\documentclass[preprint,number,12pt]{elsarticle}
\usepackage{amsmath}
\usepackage{amssymb}
\usepackage{amsthm}
\usepackage{nicefrac}
\usepackage{paralist}
\usepackage{msp_schreiber}
\begin{document}
\begin{frontmatter}
\title{Multivariate Stirling Polynomials\\of the First and Second Kind}
\author{Alfred Schreiber}
\address{Department of Mathematics and Mathematical Education, University of Flensburg,\\Auf dem Campus 1, D-24943 Flensburg, Germany}
%\ead{alfred.schreiber@uni-flensburg.de}
\ead{info@alfred-schreiber.de}
%%
%%-------------------------
%% abstract.tex
%% Last update 2013-09-07
%% Revised
%%-------------------------
%%
\begin{abstract}
Two doubly indexed families of homogeneous and isobaric polynomials in several indeterminates are considered: the (partial) exponential Bell polynomials $\bell{n}{k}$ and a new family $\spol{n}{k}\in\ints[\row{X}{1}{n-k+1}{,}]$ such that $X_1^{-(2n-1)}\spol{n}{k}$ and $\bell{n}{k}$ obey an inversion law which generalizes that of the Stirling numbers of the first and second kind. Both polynomial families appear as Lie coefficients in expansions of certain derivatives of higher order. Substituting $D^j(\varphi)$ (the $j$-th derivative of a fixed function $\varphi$) in place of the indeterminates $X_j$ shows that both $\spol{n}{k}$ and $\bell{n}{k}$ are differential polynomials depending on $\varphi$ and on its inverse $\inv{\varphi}$, respectively. Some new light is shed thereby on Comtet's solution of the Lagrange inversion problem in terms of the Bell polynomials. According to Haiman and Schmitt that solution is essentially the antipode on the Fa\`{a} di Bruno Hopf algebra. It can be represented by $X_1^{-(2n-1)}\spol{n}{1}$. Moreover, a general expansion formula that holds for the whole family $\spol{n}{k}$ ($1\leq k\leq n$) is established together with a closed expression for the coefficients of $\spol{n}{k}$. Several important properties of the Stirling numbers are demonstrated to be special cases of relations between the corresponding polynomials. As a non-trivial example, a Schl\"omilch-type formula is derived expressing $\spol{n}{k}$ in terms of the Bell polynomials $\bell{n}{k}$, and vice versa.
\end{abstract}
%%
%% End of Abstract
%%
\end{frontmatter}
%%
%% Main text:
%%-----------
%% Section 1: Introduction
%% Last update: 2021-01-26
%%
\section{Introduction}
\subsection{Background and problem}
It is well-known that a close connection exists between iterated differentiation and Stirling numbers (see, e.\,g., \cite{jord1939,rior1958,tosc1949}). Let $s_1(n,k)$ denote the signed Stirling numbers of the first kind, $s_2(n,k)$ the Stirling numbers of the second kind, and $D$ the operator $d/dx$. Then, for all positive integers $n$, the $n$th iterate $(x D)^n$ can be expanded into the sum
\begin{equation}\label{equation1.1}
	(x D)^n = \sum_{k=1}^{n} s_2(n,k) x^k D^k.
\end{equation}
An expansion in the reverse direction is also known to be valid (see, e.\,g., \cite[p.\,197]{jord1939} or \cite[p.\,45]{rior1958}):
\begin{equation}\label{equation1.2}
	D^n = x^{-n} \sum_{k=1}^{n} s_1(n,k) (x D)^k.		
\end{equation}
Let us first look at \eq{equation1.1}. The occurrence of the Stirling numbers can be explained combinatorially as follows. Observing 
\[
(x D)^n f(x)=D^n(f \circ \exp)(\log x)
\]
we can use the classical higher-order chain rule (named after Fa\`{a} di Bruno; cf. \cite{john2002,jord1939}, \cite[pp.\,52,\,481]{knut1997}) to calculate the $n$th derivative of the composite function $f\circ g$:
\begin{equation}\label{equation1.3}
	(f \circ g)^{(n)}(x)= \sum_{k=1}^n \bpol{n}{k}(g'(x),\ldots,g^{(n-k+1)}(x))\cdot f^{(k)}(g(x)),
\end{equation}
where $\bpol{n}{k}\in\ints[\row{X}{1}{n-k+1}{,}]$, $1\leq k \leq n$, is the (partial) exponential Bell polynomial
\begin{equation}\label{equation1.4}
	\bpol{n}{k}(X_{1},\ldots,X_{n-k+1}) = \sum_{r_1,r_2,\ldots}\frac{n!}{r_1!r_2!\ldots (1!)^{r_1}(2!)^{r_2}\ldots}	X_1^{r_1}X_2^{r_2}\ldots
\end{equation}
the sum to be taken over all non-negative integers $r_1,r_2,\ldots,r_{n-k+1}$ such that $r_1+r_2+\ldots+r_{n-k+1}=k$ and $r_1+2r_2+\ldots+(n-k+1)r_{n-k+1}=n$. The coefficient in $\bpol{n}{k}$ counts the partitions of $n$ distinct objects into $k$ blocks (subsets) with $r_j$ blocks containing exactly $j$ objects ($1\leq j\leq n-k+1$). Therefore, the sum of these coefficients is equal to the number $s_2(n,k)$ of all such partitions. So we have $\bpol{n}{k}(x,\ldots,x)=s_2(n,k)x^k$. Evaluating $(f \circ \exp)^{(n)}(\log x)$ by \eq{equation1.3} then immediately gives the right-hand side of \eq{equation1.1}.

\textbf{Question.} Can also \eq{equation1.2} be interpreted in this way by substituting $j$th derivatives in place of the indeterminates $X_j$ of some polynomial $\spol{n}{k}\in\ints[\row{X}{1}{n-k+1}{,}]$, the coefficients of which add up to $s_1(n,k)$?
\vspace*{1ex}

\fussy 
%\noindent
The main purpose of the present paper is to give a positive and comprehensive answer to this question including recurrences, a detailed study of the inverse relationship between the polynomial families $\bpol{n}{k}$ and $\spol{n}{k}$, as well as fully explicit formulas (with some applications to Stirling numbers and Lagrange inversion).

The issue turns out to be closely related to the problem of generali\-zing \eq{equation1.1}, that is, finding an expansion for the operator $(\theta D)^n$ ($n\geq1$, $\theta$ a function of $x$). Note that, in the case of scalar functions, $(\theta D)f$ is the \textit{Lie derivative} of $f$ with respect to $\theta$. Several authors have dealt with this problem. In \cite{comt1973} and \cite{mima1998} a polynomial family $C_{n,k}\in\ints[X_0,X_1,\ldots,X_{n-k}]$ has been defined%%
\footnote{Here we write $C_{n,k}$ instead of Comtet's $A_{n,k}$ (cf. \cite{comt1973}) in order to avoid misunderstandings. Note that in the present paper $A_{n,k}$ is exclusively used  to denote the `Lie coefficients' according to Todorov (see \eq{equation1.5} below).} 
by differential recurrences and shown to comply with $(\theta D)^n=\sum_{k=1}^{n}C_{n,k}(\theta,\theta',\ldots,\theta^{(n-k)})D^k$. Comtet \cite{comt1973} has tabulated $C_{n,k}$ up to $n=7$ and stated that $C_{n,k}(x,\ldots,x)=c(n,k)x^n$, where $c(n,k):=|s_1(n,k)|$ denotes the signless Stirling numbers of the first kind (`cycle numbers' according to the terminology in \cite{knut1992b}). Since however all coefficients of $C_{n,k}$ are positive, $C_{n,k}$ does not appear to be a suitable companion for $\bpol{n}{k}$ with regard to the desired inversion law.

\fussy
Todorov \cite{todo1981,todo1985} has studied the above Lie derivation with respect to a function $\theta$ of the special form $\theta(x)=\nicefrac{1}{\varphi'(x)}$, $\varphi'(x)\neq 0$. His main results in \cite{todo1981} ensure the existence of $\spol{n}{k}\in\ints[\row{X}{1}{n-k+1}{,}]$ such that 
\begin{equation}\label{equation1.5}
	\left(\varphi'(x)^{-1}D\right)^n \negthickspace f(x) = \sum_{k=1}^{n}\sfun{n}{k}(\varphi'(x),\ldots,\varphi^{(n-k+1)}(x))\cdot f^{(k)}(x),
\end{equation}
\noindent
where $\sfun{n}{k}:=X_1^{-(2n-1)}\spol{n}{k}$. While differential recurrences for $\sfun{n}{k}$ can readily be derived from \eq{equation1.5} (cf. \cite[Equation (27)]{todo1981} or a slightly modified version in \cite[Theorem~2]{todo1985}), a simple representation for $\spol{n}{k}$\,---\,as is \eq{equation1.4} for $\bpol{n}{k}$\,---\,was still lacking up to now. Todorov \cite[p.\,224]{todo1981} erroneously believed that the somewhat cumbersome `explicit' expression in \cite{comt1973} for the coefficients of $C_{n,k}$ would directly yield the coefficients of $\spol{n}{k}$. Also the determinantal form presented in \cite[Theorem~6]{todo1981} for $(D/\varphi')^n$ (and thus also for $\spol{n}{k}$) may only in a modest sense be regarded as explicit.

Nevertheless, Todorov's choice ($\theta = \nicefrac{1}{\varphi'}$) eventually proves to be a crucial idea. Among other things, it reveals that $\sfun{n}{k}$ (and thus $\spol{n}{k}$) is connected with the classical Lagrange problem of computing the compositional inverse $\inv{f}$ of a given series $f(x)=\sum_{n\geq1}(f_{n}/n!)x^n$, $f_{1}\neq0$. As we shall see later, the Taylor coefficients $\inv{f}_n$ of $\inv{f}(x)$ can be expressed simply by applying $A_{n,1}$ to the coefficients of $f$ as follows:
\begin{equation}\label{equation1.6}
	\inv{f}_n = \sfun{n}{1}(f_1,\ldots,f_n).
\end{equation}
On the other hand, Comtet \cite{comt1974} found an inversion formula that expresses $\inv{f}_n$ in terms of the (partial) exponential Bell polynomials:
\begin{equation}\label{equation1.7}
	\inv{f}_n = \sum_{k=0}^{n-1}(-1)^k f_1^{-n-k}\bpol{n+k-1}{k}(0,f_2,\ldots,f_n).
\end{equation}
This result has been shown by Haiman and Schmitt \cite{hasc1989,schm1994} to provide essentially both a combinatorial representation and a cancellation-free computation of the antipode on a Fa\`{a} di Bruno Hopf algebra (a topic that has received a lot of attention in quantum field theory due to its application to renormalization; cf. \cite{krei1998,cokr1998,figr2005}). Combining \eq{equation1.6} with \eq{equation1.7} we obtain an expression for $\sfun{n}{1}$ in terms of the Bell polynomials. This suggests looking for a similar representation for the whole family $\sfun{n}{k}$, $1\leq k \leq n$. As a main result (Theorem 6.1), we shall prove the following substantially extended version of \eq{equation1.6} \& \eq{equation1.7}: 
\begin{equation}\label{equation1.8}
	\sfun{n}{k} = \sum_{r=k-1}^{n-1}(-1)^{n-1-r}\binom{2n-2-r}{k-1}X_1^{-(2n-1)+r}\comt{2n-1-k-r}{n-1-r}.
\end{equation}
The tilde over $B$ indicates that $X_1$ has been replaced by $0$. From \eq{equation1.8} we eventually get the desired explicit standard representation for $A_{n,k}$ that corresponds to the one for $\bpol{n}{k}$ given in \eq{equation1.4}. 

%%\noindent 
Equation \eq{equation1.8} states a somewhat intricate relationship between the families $\sfun{n}{k}$ and $\bpol{n}{k}$. A simpler connection of both expressions is the following basic inversion law which generalizes the orthogonality of the Stirling numbers (cf. Section~5):
\begin{equation}\label{equation1.9}
	\sum_{j=k}^{n}\sfun{n}{j}\bpol{j}{k}=\kron{n}{k} \quad (1\leq k \leq n),
\end{equation}
%%
%\sloppy
where $\kron{n}{n}=1$, $\kron{n}{k}=0$ if $n\neq k$ (Kronecker symbol). 
\subsection{Terminology and notation}
Considering \eq{equation1.9} and the fact that the sum of the coefficients of $\sfun{n}{k}$ and of $\bpol{n}{k}$ are equal to $s_1(n,k)$ and to $s_2(n,k)$, respectively, it may be justified to call $\sfun{n}{k}$ and  $\bpol{n}{k}$ \emph{multivariate Stirling polynomials of the first and second kind}. There should be no risk of confusing them with polynomials \emph{in one variable} like those introduced and named after Stirling by Nielsen \cite{niel1904,niel1905}, neither with the closely related `Stirling polynomials' \mbox{$f_k(n):=s_2(n+k,n)$} and $g_k(n):=c(n,n-k)$ Gessel and Stanley \cite{gest1978} have investigated as functions of $n\in\ints$. 

A sequence $r_1,r_2,r_3,\ldots$ of non-negative integers is said to be an $(n,k)$-\textit{partition type}, $0\leq k\leq n$, if $r_1+r_2+r_3+\ldots=k$ and $r_1+2r_2+3r_3+\ldots=n$. The set of all $(n,k)$-\pt s is denoted by $\ptsi{n}{k}$; we write $\ptsa$ for the union of all $\ptsi{n}{k}$. In the degenerate case ($k=0$) set $\ptsi{n}{0}=\emptyset$, if $n>0$, and $\ptsi{0}{0}=\{0\}$ otherwise. Let $k\geq 1$. Since $n-k+1$ is the greatest $j$ such that $r_j>0$, we often write $(n,k)$-partition types as ordered $(n-k+1)$-tuples $(r_1,\ldots,r_{n-k+1})$. 

The polynomials to be considered in the sequel have the form 
\[
	P_{\pi}=\sum \pi(r_1,r_2,\ldots)X_1^{r_1}X_2^{r_2}\ldots, 
\]	
where the sum ranges over all elements $(r_1,r_2,\ldots)$ of a full set $\ptsi{n}{k}$. As a consequence, $P_{\pi}$ is homogeneous of degree $k$ and isobaric of degree $n$. The coefficients of $P_{\pi}$ may be regarded as values of a map $\pi:\ptsa\longrightarrow\ints$ defined by some combinatorially meaningful expression, at least in typical cases like the following:
{\allowdisplaybreaks
\begin{alignat}{6}
 \ord(r_1,r_2,\ldots):=& \dfrac{(r_1+2r_2+\ldots)!}{r_1!\,r_2!\cdot\ldots~~~}
      & \text{\textit{order function} (Lah)}\label{lah}\\[6pt]
 \cau(r_1,r_2,\ldots):=& \dfrac{\omega(r_1,r_2,\ldots)}{1^{r_1}\,2^{r_2}\cdot\ldots~}   
      & \text{\textit{cycle function} (Cauchy)}\label{cauchy}\\[6pt]
 \bru(r_1,r_2,\ldots):=& \dfrac{\omega(r_1,r_2,\ldots)}{(1!)^{r_1}\,(2!)^{r_2}\cdot\ldots}   
      & \text{\textit{~~subset function} (Fa\`{a} di Bruno)}\label{bruno}
\end{alignat}
}% end of \allowdisplaybreaks
These coefficients count the number of ways a set can be partitioned into non-empty blocks according to a given \pt, that is, $r_{j}$ denotes the number of blocks containing $j$ elements ($j=1,2,\ldots$). The result depends on the meaning of `block': linearly ordered subset ($\ord$), cyclic order ($\cau$), or unordered subset ($\bru$). 

\sloppy
It should be noticed that the corresponding polynomials $P_{\ord},P_{\cau},P_{\bru} (=\bpol{n}{k})$ are closely related to well-known combinatorial number-families: 
\begin{align*}
	&P_{\ord}(1,\ldots,1)=l^{+}(n,k),\text{~~unsigned Lah numbers \cite{lah1955,rior1958}}\\
	&P_{\cau}(1,\ldots,1)=c(n,k)={\cycs{n}{k}},\text{~~unsigned Stirling numbers of the 1st kind}\\
	&P_{\bru}(1,\ldots,1)=s_2(n,k)={\subs{n}{k}},\text{~~Stirling numbers of the 2nd kind}. 
\end{align*}

\fussy
\subsection{Overview} 
This paper is organized as follows: In Section~2 a general setting is sketched that allows functions and derivations to be treated algebraically. Section~3 contains a study of the iterated Lie operator $D(\varphi)^{-1}D$. An expansion formula for $(D(\varphi)^{-1}D)^n$ is established together with a differential recurrence for the resulting Lie coefficients $\sfun{n}{k}$. Doing the same with respect to the inverse function $\inv{\varphi}$ will yield, conversely, $D^n$ expanded and $\bpol{n}{k}$ as the corresponding Lie coefficients. A by-product of Section~3 is Fa\`{a} di Bruno's formula and its applications to the partial Bell polynomials $\bpol{n}{k}$ to be briefly  summarized in Section~4. These basic facts then lead to both inversion and recurrence relations, which we shall demonstrate and discuss in Section~5. The main task in Section~6 is to find an explicit polynomial expression for $\spol{n}{k}$. This is eventually achieved by means of \eq{equation1.8}, a proof of which makes up a central part of the section. In Section~7 we give some applications to the Lagrange inversion problem and to exponential generating functions.
%%
%% End of text
%%%%%%%%%%%%%%%%%%%%%%%%%%%% 
%% Section 2: Function algebra with derivation
%% Last updated: 2021-01-26
%
\section{Function algebra with derivation}
\subsection{Basic notions}
Menger \cite{meng1944} has introduced the notion of a `tri-operational algebra' of functions, which in the sequel (since 1960) stimulated to a great extent studies of generalized function algebras, e.\,g., \cite{dutr2012,meng1961,schw1960,schw1967,suci1975}. In what follows I will propose a variant of Menger's original system tailored to our specific purposes of treating functions and their derivatives in a purely algebraic way.

Let $(\funs,+,\cdot)$ be a non-trivial commutative ring of characteristic zero, 0 and 1 its identity elements with respect to addition and multiplication. We will think of the elements of $\funs$ as `functions (of one variable)' and therefore assume that $\funs$ has a third binary operation $\circ$ (called \textit{composition}) together with an identity element $\id$ such that the following axioms are satisfied:
\vspace*{-2pt}
{\allowdisplaybreaks
\begin{align*}
\text{(F1)} \quad & f \circ (g \circ h) = (f \circ g) \circ h \mspace{280mu}\\
\text{(F2)} \quad & (f+g)\circ h = (f \circ h) + (g \circ h)\\
\text{(F3)} \quad & (f \cdot g)\circ h = (f \circ h) \cdot (g \circ h)\\
\text{(F4)} \quad & f \circ \iota = \iota \circ f = f\\
\text{(F5)} \quad & 1 \circ 0 = 1\\
\end{align*}
}% end of \allowdisplaybreaks

\vspace*{-4ex}
\noindent
(F4) is assumed to be valid for all $f\in\funs$; hence $\iota$ is unique. Let $f$ be any element of $\funs$. From (F2) we conclude $0\circ f=0$; so we get $\iota\neq0$ (by (F4)) and $\iota\neq1$ (by (F5)). (F2) furthermore implies $(-f)\circ g=-(f\circ g)$. 

The least subring of $\funs$ containing 1 will in the following conveniently be identified with $\ints$. (F5) then extends to the integers, that is, $n\circ 0=n$ holds for all $n\in\ints$.

Given a \textit{unit} $f$ in $\funs$ (i.\,e., $f$ is an element invertible with respect to multiplication), we write $f^{-1}$ (or $\nicefrac{1}{f}$) for the inverse (henceforth called \textit{reciprocal}) of $f$. 

\smallskip
\begin{rem}
	\label{rem2.1} 
	It must be emphasized that $\circ$ has to be understood as a \textit{partial} operation (of course, $\iota^{-1}\circ\,0$ is not defined). We therefore assign truth values to formulas, especially to our postulates (F1--3), whenever the terms involved are meaningful. 
\end{rem}
\smallskip

\noindent
Let $f,g\in\funs$ be functions such that $f \circ g = g \circ f = \iota$. Then $g$ is called the \textit{compositional inverse} of $f$, and vice versa. It is unique and will be denoted by $\inv{f}$. The following is obvious: $\inv{\iota}=\iota$, $\inv{\inv{f}}=f$, and $\inv{f \circ g}=\inv{g} \circ \inv{f}$.

\smallskip
\begin{dfn}
\label{dfn2.1}
Suppose $(\funs,+,\cdot,\circ)$ satisfies (F1--5). We then call a mapping $D: \funs \longrightarrow \funs$ \textit{derivation on} $\funs$, and $(\funs,+,\cdot,\circ,D)$ a \textit{function algebra with derivation}, if $D$ meets the following conditions:
{\allowdisplaybreaks
\begin{align*}
\text{(D1)} \quad & D(f+g)=D(f)+D(g) \mspace{280mu}\\
\text{(D2)} \quad & D(f \cdot g)=D(f)\cdot g + f\cdot D(g)\\
\text{(D3)} \quad & D(f \circ g) = (D(f) \circ g)\cdot D(g)\\
\text{(D4)} \quad & D(\iota)=1\\
\text{(D5)} \quad & D(f)=0\Longrightarrow f\circ 0=f
\end{align*}
}% end of \allowdisplaybreaks
\end{dfn}
\label{Product rule}

\noindent
The classical derivation rules (D1), (D2) make $\funs$ into a differential ring. Some simple facts are immediate: $D(0)=D(1)=0$, $D(m\cdot f)=m\cdot D(f)$ for all $m\in\ints$. By an inductive argument the product rule (D2) can be generalized:
\begin{equation}\label{equation2.1}
	D(f_1\dotsm f_n)=\sum_{k=1}^{n}f_1\cdots f_{k-1}{\,}\cdot D(f_k)\cdot f_{k+1}\cdots f_n.
\end{equation}
\label{Logarithmic derivative identity}

\noindent
By putting $f_i=f$, $1\leq i\leq n$, \eq{equation2.1} becomes $D(f^n)=n f^{n-1} D(f)$.  If $f$ is a unit, this holds also for $n\leq0$. As usual, $f^m$  for $m<0$ is defined by $(f^{-1})^{-m}$.

(D4) prevents $D$ from operating trivially. In the case of a field $\funs$, (D4) can be weakend to $D(f)\neq0$ (for some $f\in\funs$), since the chain rule (D3) then gives $D(f)=D(f\circ\iota)=D(f)\cdot D(\iota)$. 

Applying (D3) and (D4) to $D(f\circ\inv{f})$ we obtain the inversion rule
\begin{equation}\label{equation2.2}
	D( \inv{f}) = \frac{1}{D(f) \circ \inv{f}}.
\end{equation}
\label{Notion of constant}

\noindent\sloppy
In a differential ring, it is customary to define the subring $\cons$ of constants as the kernel of the additive homomorphism $D$, that is, 
\[
	\cons:=\{f\in\funs\,|\,D(f)=0\}. 
\]	
We have $\ints\subseteq\cons$. Constants behave as one would expect.

\smallskip\fussy
\begin{prop}
  \label{prop2.1}
  $c\in\cons\iff c\circ f=c$~~for all $f\in\funs$.
\end{prop}	
\begin{proof}
$\Rightarrow$: Suppose $D(c)=0$. Then, for any $f\in\funs$ we have by (F1) and (D5): $c\circ f=(c\circ 0)\circ f=c\circ(0\circ f)=c\circ 0=c$. --- $\Leftarrow$: Set $f=0$ and apply the chain rule (D3).
\end{proof}

If $f\circ 0$ exists for $f\in\funs$, then it is obviously a constant.

\sloppy
\begin{exms}[Function algebras with derivation]\label{ex2.1}
In each case below, the `functions' carry some argument $X$ (indeterminate, variable) that can be substituted in the usual sense: $f(g(X))=:(f \circ g)(X)$. 
\begin{compactenumerate}
	\item The rational function field $\rats(X)$ together with the algebraically defined derivation $R \longmapsto R'$, $R\in\rats(X)$.
	\item The ring $\real[[X]]$ of all power series with formal differentiation.
  \item The ring of real-valued $C^{\infty}$ functions on an open real interval with the ordinary differential operator $d/dx$.
	\item The field of all meromorphic functions on a given region in $\complex$ with complex differentiation.
\end{compactenumerate}
\end{exms}

\fussy
%\noindent
Given a polynomial $P\in\cons[\row{X}{1}{n}{,\,}]$ and functions $\row{f}{1}{n}{,\,}$, we denote by $P(\row{f}{1}{n}{,\,})$ the function obtained by substituting $f_i$ in place of $X_i$, $1\leq i\leq n$. Recalling the algebraic definition of $\partial/\partial X_i$, we readily obtain by \eq{equation2.1} the generalized chain rule:
\begin{equation}\label{equation2.3}
	D(P(\row{f}{1}{n}{,})) = \sum_{k=1}^n \pder{P}{X_k}\,(\row{f}{1}{n}{,\,})\cdot D(f_k).
\end{equation}
\label{Generalized chain rule}

\noindent\fussy 
In the case $n=1$ this becomes $D(P(f))=P'(f)\cdot D(f)$. Thus, $D$ restricted to polynomial functions turns out to act like the ordinary differential operator: $D(P(\iota))=P'(\iota)$. Equation \eqref{equation2.3} also applies to the case that $P$ is a rational function from $\cons(X_1,\ldots,X_n)$.
\begin{ntn}
Suppose that $\varphi$ is a fixed function, and let $D^i$ denote the $i$th iterate of $D$. In the following we will abbreviate $P(D(\varphi),D^2(\varphi),\,\ldots\,,D^n(\varphi))$ to $P^{\dph}$. We denote by $\poln_n(\varphi)$ the set of all $P^{\dph}$ and by $\ratf_n(\varphi)$ the set of all rational expressions $P^{\dph}/Q^{\dph}$, where $P,Q\in\cons[X_1,\ldots,X_n]$.
\end{ntn}
\begin{rem}\label{rem2.3}
Obviously $(P+Q)^{\dph}=P^{\dph}+Q^{\dph}$ and $(P\cdot Q)^{\dph}=P^{\dph}\cdot Q^{\dph}$. For a homogeneous polynomial $P$ of degree $k$ and arbitrary $a,b\in\cons$ we have $P^{\mspace{2mu}a\varphi+b}=a^k P^{\dph}$.
\end{rem}
\label{Homogeneous polynomials}
\sloppy
\noindent
Another formula that is known from calculus and that holds in the differential ring $(\funs,D)$ is the general Leibniz rule, which yields an explicit expression for the higher derivatives of a product:
\begin{equation}\label{equation2.4}
	D^s(f_1\cdots f_n)=\sum_{\substack{j_1+\cdots+j_n=s \\ 
			                                        j_1,\,\ldots,j_n\,\geq\,0}}\,
	\frac{s!}{j_1!\cdots j_n!}\;D^{j_1}(f_1)\cdots D^{j_n}(f_n).
\end{equation}

\fussy
We note the special case $D^s(\iota^n)=\powfall{n}{s}\cdot\iota^{n-s}$, if $s\leq n$, and $D^s(\iota^n)=0$ otherwise; the falling power $\powfall{n}{s}$ is defined by $n(n-1)\cdots(n-s+1)$.
\begin{cnv}
\textit{Throughout the remainder of this paper we denote by $\varphi$ any function from $\funs$ such that the compositional inverse $\inv{\varphi}$ exists, $D(\varphi)$ and $D(\inv{\varphi})$ are units, and equation \eq{equation2.2} holds for $f=\varphi$.}
\end{cnv}
With regard to such a $\varphi$, we define a mapping $D_{\varphi}:\funs\longrightarrow\funs$ by
\begin{equation}\label{equation2.5}
	\lieder{\varphi}{f}:=\frac{D(f)}{D(\varphi)}.
\end{equation}

\noindent
$D_{\varphi}$ is the function-algebraic version of the \textit{Lie derivative} that Todorov used in his paper \cite{todo1981}. $(\funs,D_{\varphi})$ is a differential ring having the same constants as $(\funs,D)$. However, with regard to $D_{\varphi}$ we have: (D3)$\;\Leftrightarrow\;$(D4)$\;\Leftrightarrow D(\varphi)=1$. Therefore, $D_{\varphi}$ satisfies the chain rule if and only if $D_{\varphi}=D$, that is, in the \textit{trivial} case $\varphi=\iota + c$, $c\in\cons$.

The following simple but useful statement concerns the relation between $D_{\varphi}^n$ and $D^n$ (\textit{Pourchet's formula} according to \cite[p.\,220]{comt1974} and \cite[p.\,223--224]{todo1981}).
\begin{prop}\label{prop2.2}
	$D_{\varphi}^n(f) = D^n(f \circ \inv{\varphi}) \circ \varphi \qquad \text{for all~}n\geq 0.$
\end{prop}
\begin{proof}
Verify $D^n(f\circ\inv{\varphi})=D_{\varphi}^n(f)\circ\inv{\varphi}$ by induction on $n$. The case $n=0$ is clear. For the induction step ($n\rightarrow n+1$) use (D3): 
\[
	D^{n+1}(f\circ\inv{\varphi})=(D(D_{\varphi}^n(f))\circ\inv{\varphi}))\cdot D(\inv{\varphi})=:(*). 
\]
Applying \eq{equation2.2} to $D(\inv{\varphi})$ yields $(*)=\frac{D(D_{\varphi}^n(f))}{D(\varphi)}\circ\inv{\varphi}=D_{\varphi}^{n+1}(f)\circ\inv{\varphi}$.\end{proof}

\subsection{Exponential and logarithm}
For some purposes it will prove convenient to have in $(\funs,+,\cdot,\circ,D)$ besides the three identity elements more functions with special properties, the most important examples being the exponential ($\exp$) and its compositional inverse ($\log$). Of course we expect $\exp$ and $\log$ to have the familiar pro\-perties known from analysis, like $D(\exp)=\exp$, $D(\log)\cdot\id=1$, and $D^k(\log)=(-1)^{k-1}(k-1)!\,\id^{-k}$ for all $k\geq 1$. In such a case we call $\funs$ an \emph{extended function algebra}. Items 3 and 4 from \examples{ex2.1} are extended function algebras. In $\real[[x]]$ (item~2), however, neither $D(\log)$ nor $\log$ have counterparts. In Section 7 we shall therefore deal with $\log\circ(1+\iota)$ and $\exp-1$ instead. 

\begin{ntn}
In an extended function algebra, we prefer to write $\idty$ (or $x$) instead of $\id$. If no misunderstanding is likely, we will occasionally replace $\exp\circ f$ with $e^f$ and also switch from $f\circ g$ to the usual notation $f(g)$, mainly when $ g $ is a constant or the identity function.
\end{ntn}
%
% End of text
%%%%%%%%%%%%%%%%%%%%%%%%%%%%
%% Section 3: Expansion of higher order derivatives
%% Last updated: 2020-12-05
%
\section{Expansion of higher-order derivatives}
\sloppy
\noindent
In this section, some basic facts regarding the expansion of $D_{\varphi}^n$ (see \eq{equation1.5} and Todorov \cite{todo1981,todo1985}) as well as of $D^n$ will be reformulated and set up for an arbitrary function algebra $\funs$ with derivation $D$. The main idea of the presentation is to make it clear that from the beginning these results are linked by an inversion of functions. We will also examine some basic properties of the multivariable polynomials involved, in particular their differential recurrences and their relationship to the Stirling numbers.
\subsection{Expansion formulas for $D_{\varphi}^{n}$ and $D^{n}$}

\begin{prop}\label{prop3.1}
	Let $f$ be any function from $\funs$ and $n,k$ non-negative integers. Then there are $a_{n,k}\in\ratf_{n-k+1}(\varphi)$, $0\leq k\leq n$, such that
\[
	\hlieder{\varphi}{n}{f} = \sum_{k=0}^{n} a_{n,k}\cdot D^{k}(f). 
\]
The coefficients $a_{n,k}$ are uniquely determined by the recurrence
\[
	\quad a_{n+1,k}=\frac{a_{n,k-1}+D(a_{n,k})}{D(\varphi)} \quad (1\leq k\leq n+1),
\]
where $a_{0,0}=1$, $a_{i,0}=0$ $(i>0)$, and $a_{i,j}=0$ $(0\leq i<j)$.
\end{prop}
\fussy
\begin{proof}
Recall that $D(\varphi)$ is a unit and $(\funs,D_{\varphi})$ is a differential ring. So, $\hlieder{\varphi}{n}{f}$ is defined for every $n\geq 0$ and can successively be calculated by applying (D1) and (D2) together with the rule $D_{\varphi}(g^{-1})=-g^{-2}D_{\varphi}(g)$ (\mbox{$g$ a unit}). The proof is then carried out by a simple induction on $n$, the details of which can be omitted here.
\end{proof}
\begin{rem}\label{rem3.1}
$a_{n,n}=D(\varphi)^{-n}$ for all $n \geq 0$.
\end{rem}
\noindent
One obtains by induction that the denominator in $a_{n,k}$ is $D(\varphi)^{2n-1}$. 
\sloppy
\begin{cor}\label{cor3.2}
Set $s_{n,k}:=D(\varphi)^{2n-1}a_{n,k}$. Then $s_{n,k}\in\poln_{n-k+1}(\varphi)$ for \mbox{$(n,k)\neq (0,0)$}, and the following recurrence holds: 
\[
s_{n+1,k}=-(2n-1)D^2(\varphi)s_{n,k}+D(\varphi)\cdot(s_{n,k-1}+D(s_{n,k})),
\]
where $s_{0,0}=D(\varphi)^{-1}$, $s_{1,1}=1$, $s_{i,0}=0$ $(i>0)$, and $s_{i,j}=0$ $(0\leq i<j)$.
\end{cor}

\fussy
It is natural to ask whether, conversely, $D^{n}$ can be expanded into a linear combination of the $D_{\varphi}^{k}$ $(k=0,1,\ldots,n)$. The next proposition gives a positive answer.
\begin{prop}\label{prop3.3}
	Let $f$ be any function and $n,k$ non-negative integers. Then, there are $b_{n,k}\in\poln_{n-k+1}(\varphi)$, $0\leq k\leq n$, such that
\[
	D^{n}(f) = \sum_{k=0}^{n} b_{n,k}\cdot \hlieder{\varphi}{k}{f}. 
\]
The coefficients $b_{n,k}$ are uniquely determined by the recurrence
\[
	\quad b_{n+1,k}={D(\varphi)}\cdot b_{n,k-1}+D(b_{n,k}) \quad (1\leq k\leq n+1),
\]
where $b_{0,0}=1$, $b_{i,0}=0$ $(i>0)$, and $b_{i,j}=0$ $(0\leq i<j)$.
\end{prop}
\begin{proof}
\sloppy
We apply \proposition{prop3.1} to the compositional inverse $\inv{\varphi}$, thus obtaining $a'_{n,k}\in\ratf_{n-k+1}(\inv{\varphi})$ so that $\hlieder{\inv{\varphi}}{n}{f} = \sum_{k=0}^{n} a'_{n,k}\cdot D^{k}(f)$. Since according to Pourchet's formula (\proposition{prop2.2}) the left-hand side is equal to \mbox{$D^n(f\circ\varphi)\circ\inv{\varphi}$}, we get by (F3)
\begin{equation}\label{equation3.1}
	D^n(f\circ\varphi) = \sum_{k=0}^{n}(a'_{n,k}\circ\varphi)\cdot(D^{k}(f)\circ\varphi).
\end{equation}
\fussy
We now replace $f$ by $f\circ\inv{\varphi}$ in \eq{equation3.1}. Then again Pourchet's formula, applied to the second factor on the right-hand side of \eq{equation3.1}, yields 
\[
	D^{n}(f)=\sum_{k=0}^{n}(a'_{n,k}\circ\varphi)\cdot \hlieder{\varphi}{k}{f}.
\]
Now set $b_{n,k}:=a'_{n,k}\circ\varphi$. We then have $b_{0,0}=a'_{0,0}\circ\varphi=1\circ\varphi=1$, likewise $b_{i,0}=0$ $(i>0)$, $b_{i,j}=0$ $(0\leq i<j)$, and by \eq{equation2.2} 
{\allowdisplaybreaks
\begin{align*}
	b_{n+1,k} &= a'_{n+1,k}\circ\varphi=(D(\inv{\varphi})^{-1}\circ\varphi)\cdot((D(a'_{n,k})\circ\varphi)+(a'_{n,k-1}\circ\varphi))\\
	          &= D(\varphi)\cdot (D(a'_{n,k})\circ\varphi) + D(\varphi)\cdot b_{n,k-1}\\
						&= D(b_{n,k}) + D(\varphi)\cdot b_{n,k-1}\qquad\qquad~~~~~~~~\text{(by the chain rule (D3))}.
\end{align*}
}% Ende allowdisplaybreaks
Finally, $b_{n,k}\in\poln_{n-k+1}(\varphi)$ follows from this recurrence by an inductive argument.
\end{proof}
\begin{rem}\label{rem3.2}
$b_{n,n}=D(\varphi)^{n}$ for all $n \geq 0$.
\end{rem}

\subsection{Fundamental properties of the coefficients}
\sloppy
Let us now have a closer look at the coefficient functions $a_{n,k}$ and $b_{n,k}$. We start with $b_{n,k}$. As it is a polynomial expression in the derivatives $D(\varphi),\ldots,D^{n-k+1}(\varphi)$, we get $b_{n,k}$ from a suitable polynomial family $\bpol{n}{k}$ by substituting $D^{j}(\varphi)$ in place of the indeterminates $X_{j}$, that is, $\bpol{n}{k}^{\dph}=b_{n,k}$. In the case of $a_{n,k}$ it is likewise clear by \corollary{cor3.2} that $s_{n,k}$, too, comes from certain polynomials $\spol{n}{k}$ satisfying $\spol{n}{k}^{\dph}=s_{n,k}=D(\varphi)^{2n-1}a_{n,k}$. We therefore define $\sfun{n}{k}:=X_{1}^{-(2n-1)}\spol{n}{k}$. Then, of course, $\sfun{n}{k}^{\dph}=a_{n,k}$ holds. Note that $A_{n,k}$ is a Laurent polynomial, and that is especially also true for $\spol{0}{0}$ (see \corollary{cor3.2}).

\fussy
The polynomials $\spol{n}{k}$ ($A_{n,k}$) and $\bpol{n}{k}$ are closely connected.

\begin{prop}\label{prop3.4}
{\allowdisplaybreaks
\begin{alignat*}{2}
	\text{\textup{(i)}}\qquad
	   \bpol{n}{k}^{\dph} &= D(\varphi)^{2n-1}\cdot(\spol{n}{k}^{\dphinv}\circ\varphi),\qquad & \bpol{n}{k}^{\dph} = \sfun{n}{k}^{\dphinv}\circ\varphi\\
   \text{\textup{(ii)}}\qquad
	   \spol{n}{k}^{\dph} &= D(\varphi)^{2n-1}\cdot(\bpol{n}{k}^{\dphinv}\circ\varphi),\qquad & \sfun{n}{k}^{\dph} = \bpol{n}{k}^{\dphinv}\circ\varphi
\end{alignat*}
}% end of \allowdisplaybreaks
\end{prop}
\begin{proof}
(i): From the proof of \proposition{prop3.3} we obtain 
\[
\bpol{n}{k}^{\dph}=a'_{n,k}\circ\varphi = \sfun{n}{k}^{\dphinv}\circ\varphi. 
\]
Now note the generalized inversion rule obtained by induction from \eq{equation2.2}: 
\[
	D(\inv{\varphi})^{m}\circ\varphi=D(\varphi)^{-m} \quad\text{for all integers~}m\geq 0.
\]
\label{Generalized inversion rule}

\vspace*{-1.7ex}\noindent
Hence $\spol{n}{k}^{\dphinv}\circ\varphi=(D(\inv{\varphi})^{2n-1}\circ\varphi)\cdot(\sfun{n}{k}^{\dphinv}\circ\varphi)=D(\varphi)^{-(2n-1)}\cdot\bpol{n}{k}^{\dph}$.\\(ii): Replace $\varphi$ in (i) by $\inv{\varphi}$.
\end{proof}

%%\noindent 
From the foregoing we gather the following special values: 
{\allowdisplaybreaks
\begin{align*}
  &\sfun{n}{n} = X_{1}^{-n}, \quad \spol{n+1}{n+1}=\bpol{n}{n}=X_{1}^{n}\quad &(n\geq 0),\\
	&\sfun{i}{0} =\spol{i}{0}=\bpol{i}{0}=0\quad &(i>0),\\
	&\sfun{i}{j} =\spol{i}{j}=\bpol{i}{j}=0\quad &(0\leq i<j).
\end{align*}
}

\vspace*{-3ex}
%%\noindent 
What still remains to be done is transforming the differential recurrences for $a_{n,k}$, $s_{n,k}$, $b_{n,k}$ into recurrences for the corresponding polynomials $\sfun{n}{k}$, $\spol{n}{k}$, $\bpol{n}{k}$. Consider the derivative 
\[
	D(a_{n,k})=D(\sfun{n}{k}(D(\varphi),\ldots,D^{n-k+1}(\varphi))).
\]
Applying \eq{equation2.3} to the right-hand side, we obtain
{\allowdisplaybreaks
\begin{align*}
	D(a_{n,k}) &= \sum_{j=1}^{n-k+1}\pder{\sfun{n}{k}}{X_{j}}\left(D(\varphi),\ldots,D^{n-k+1}(\varphi)\right)\cdot D^{j+1}(\varphi)\\
	           &= \left(\sum_{j=1}^{n-k+1}X_{j+1}\pder{\sfun{n}{k}}{X_{j}}\right)^{\varphi}.
\end{align*}
$D(s_{n,k})$ and $D(b_{n,k})$ resolve in the same manner. So, according to \proposition{prop3.1}, \corollary{cor3.2} and \proposition{prop3.3} we have the following\smallskip

\begin{prop}\label{prop3.5}
{\allowdisplaybreaks
\begin{alignat*}{2}
	\text{\textup{(i)}}\quad 
	\sfun{n+1}{k} &= \frac{1}{X_{1}}\left(\sfun{n}{k-1} + \sum_{j=1}^{n-k+1}X_{j+1}\pder{\sfun{n}{k}}{X_{j}}\right), \\ 
	\text{\textup{(ii)}}\quad 
	\spol{n+1}{k} &= -(2n-1)X_2\spol{n}{k} + X_1\left(\spol{n}{k-1} + \sum_{j=1}^{n-k+1}\mspace{-10mu}X_{j+1}\pder{\spol{n}{k}}{X_{j}}\right), \\
	\text{\textup{(iii)}}\quad 
	\bpol{n+1}{k} &= X_1\bpol{n}{k-1} + \sum_{j=1}^{n-k+1}X_{j+1}\pder{\bpol{n}{k}}{X_{j}}.
\end{alignat*}
}% end of \allowdisplaybreaks
\end{prop}
\label{Recurrences}

\sloppy
%\noindent
It follows (by induction) from \proposition{prop3.5} that these polynomials have integral coefficients. We denote by $\cst{n}{k}(r_1,\ldots,r_{n-k+1})$ and $\cbe{n}{k}(r_1,\ldots,r_{n-k+1})$ the coefficients of $X_1^{r_1}\cdots X_{n-k+1}^{r_{n-k+1}}$ in $\spol{n}{k}$ and in $\bpol{n}{k}$, respectively, thus obtaining

\fussy
\begin{cor}\label{cor3.6}
{\allowdisplaybreaks
\begin{alignat*}{2}
	\text{\textup{(i)}}\qquad 
	\spol{n}{k}~ &= \sum_{\ptsi{2n-1-k}{n-1}}\mspace{-32mu}\cst{n}{k}(r_1,\dots,r_{n-k+1})X_1^{r_1}\cdots X_{n-k+1}^{r_{n-k+1}},\\ 
	\text{\textup{(ii)}}\qquad 
	\bpol{n}{k}~ &= ~\quad\sum_{\ptsi{n}{k}}\mspace{8mu}\cbe{n}{k}(r_1,\dots,r_{n-k+1})X_1^{r_1}\cdots X_{n-k+1}^{r_{n-k+1}}.
\end{alignat*}
}% end of \allowdisplaybreaks
\end{cor}
\begin{proof}
(i): By induction on $n$. For $n=1$ we have the degenerate case of a $(0,0)$-\pt, $r_1=0$. Thus $\spol{1}{1}=1$ can be achieved by choosing \mbox{$\cst{1}{1}(0)=1$}. The induction step ($n\rightarrow n+1$) is carried out by examining the \pt s produced by the terms $X_2\spol{n}{k}$, $X_1\spol{n}{k-1}$, and $X_1X_{j+1}\pder{\spol{n}{k}}{X_{j}}$ in part (ii) of \proposition{prop3.5}. Each of them  makes $\sum r_i \;(=n-1)$ increase by 1 and makes $\sum i r_i\;(=2n-1-k)$ increase by 2, which gives the appropriate $(2n+1-k,n)$-\pt s for $\spol{n+1}{k}$. --- (ii): Obviously due to a similar argument.
\end{proof}
\begin{rem}\label{rem3.3}
We already know that $\spol{n}{n}=X_{1}^{n-1}$. Taking $k=n$ in \proposition{prop3.5}\,(ii) then yields $\spol{n}{n-1}=-\binom{n}{2}X_{1}^{n-2}X_{2}$ ($n\geq 2$). Todorov \cite{todo1985} has also calculated $\spol{n}{n-2}$ and $\spol{n}{n-3}$ this way.
\end{rem}
\begin{rem}\label{rem3.4}
As a consequence of \corollary{cor3.6}, $\spol{n}{k}$ is homogeneous of degree $n-1$ and isobaric of degree $2n-1-k$, while $\bpol{n}{k}$ is homogeneous of degree $k$ and isobaric of degree $n$.
\end{rem}

We now define integers $s_1(n,k):=\sfun{n}{k}(1,\ldots,1)=\spol{n}{k}(1,\ldots,1)$ and $s_2(n,k):=\bpol{n}{k}(1,\ldots,1)$, that is, $s_1(n,k)$, $s_2(n,k)$ are the sums of the coefficients of $A_{n,k}$ ($\spol{n}{k}$) and $\bpol{n}{k}$, respectively.

\begin{prop}\label{prop3.7}
Let $n,k$ be integers, $0\leq k\leq n$. Then, $s_1(n,k)$ are the signed Stirling numbers of the first kind, and $s_2(n,k)$ are the Stirling numbers of the second kind:
\[
	\text{\textup{(i)}}\quad s_1(n,k)=(-1)^{n-k}\cycs{n}{k},\quad\qquad\text{\textup{(ii)}}\quad s_2(n,k)=\subs{n}{k}.
\]
\end{prop}
\begin{proof}
(i): From the special values above we gather $s_1(0,0)=1$, $s_1(i,0)=0$ $(i>0)$, and $s_1(i,j)=0$ $(0\leq i<j)$. It suffices to show that $s_1$ satis\-fies the recurrence $s_1(n+1,k)=s_1(n,k-1)-n s_1(n,k)$, which defines the Stirling numbers of the first kind (see e.\,g. \cite[p.\,33]{rior1958} ). We use \corollary{cor3.6}\,(i). Consider first
\begin{equation*}
	\begin{split}
	X_{j+1}\,\pder{\spol{n}{k}}{X_{j}}=\,\sum_{\ptsi{2n-1-k}{n-1}}\mspace{-5mu}r_{j}\,\cst{n}{k}(r_1,\dots,r_{n-k+1})\,\times \\
	\times\,X_1^{r_1}\cdots X_{j}^{r_{j}-1}X_{j+1}^{r_{j+1}+1} \cdots X_{n-k+1}^{r_{n-k+1}}.
	\end{split}
\end{equation*}

\noindent
Replacing all indeterminates by 1 and then taking the sum from $j=1$ to $n-k+1$ yields
\[
 \sum_{j=1}^{n-k+1}\sum_{\ptsi{2n-1-k}{n-1}}\mspace{-20mu}r_{j}\,\cst{n}{k}(r_{1},\dots,r_{n-k+1})=s_1(n,k)\negmedspace\sum_{j=1}^{n-k+1}r_{j}. 
\]
Observing $r_{1}+\dots+r_{n-k+1}=n-1$ we get by \proposition{prop3.5}\,(ii)
\begin{align*}
s_1(n+1,k)&=-(2n-1)s_1(n,k) + s_1(n,k-1) + (n-1)s_1(n,k)\\
          &= s_1(n,k-1) - n s_1(n,k).
\end{align*}
(ii): The recurrence $s_2(n+1,k)=s_2(n,k-1)+k s_2(n,k)$ that defines the Stirling numbers of the second kind, can be verified by a similar argument using (iii) from \proposition{prop3.5}.
\end{proof}

\begin{exms}\label{ex3.5}
(i) We consider some special cases in an extended  function algebra:
\vspace*{-4pt}
\begin{alignat*}{2}
  \text{\textup{(1)}} &\qquad\dsfun{\exp}{n}{k}= s_{1}(n,k)\cdot\exp^{-n}&\qquad\text{\textup{(1')}}\qquad &\dbell{\log}{n}{k}=s_{1}(n,k)\cdot\idty^{-n}\\
	\text{\textup{(2)}} &\qquad\dbell{\exp}{n}{k}=s_{2}(n,k)\cdot\exp^{k}&\qquad\text{\textup{(2')}}\qquad &\dsfun{\log}{n}{k}= s_{2}(n,k)\cdot\idty^{k}.	
\end{alignat*}	

\vspace*{1ex}
\noindent
(1') and (2') immediately follow by \proposition{prop3.4} from (1) and (2), respectively. It is enough to perform the calculation for (1):
\vspace*{1ex}
\begin{align*}
	A_{n,k}^{\exp} &= A_{n,k}(D(\exp),\ldots,D^{n-k+1}(\exp))\\
		&=D(\exp)^{-(2n-1)}\cdot\,S_{n,k}(\exp,\ldots,\exp)\\
		&=\exp^{-(2n-1)}\cdot\,\exp^{n-1}\cdot\,S_{n,k}(1,\ldots,1) &\text{(\remark{rem3.4})}\\
		&= s_{1}(n,k)\cdot\exp^{-n}.&\text{\textup{(Def. $s_1$, \proposition{prop3.7})}} 
\end{align*}

\noindent
(ii) Since \proposition{prop3.7} tells us that $s_{1}(n,k)$ are in fact the signed Stirling numbers of the first kind, we choose $\varphi=\log$ in the expansion of \proposition{prop3.3}, which implies equation \eq{equation1.2} in the form
\[
	D^{n}(f)=\sum_{k=1}^{n}\dbell{\log}{n}{k}\cdot\hlieder{\log}{k}{f}=\idty^{-n}\sum_{k=1}^{n}s_{1}(n,k)(\idty\cdot D)^{k}(f).
\]
Analogously combining (2) and (2') with \proposition{prop3.1}, the reader should verify also the following expansion formula that corresponds to \eq{equation1.1}:
\begin{align*}
\qquad(\idty\cdot D)^{n}(f)&=\hlieder{\log}{n}{f}=\sum_{k=1}^{n}s_{2}(n,k)\cdot\idty^{k}\cdot D^{k}(f).
\end{align*}
\label{Normal ordering}

\vspace*{-2.7ex}\noindent
(iii) It may be of interest to show how some special Stirling numbers can be directly calculated within the machinery of an extended function algebra. Let us, for instance, compute $s_{1}(n,1)$. By (1) and \proposition{prop3.4}
\begin{equation*}
	s_{1}(n,1)=\exp^{n}\cdot\dsfun{\exp}{n}{1} = \exp^{n}\cdot(\dbell{\log}{n}{1}\circ\exp)=\exp^{n}\cdot(D^{n}(\log)\circ\exp).
\end{equation*}
Observing that $s_{1}(n,k)\in\cons$, we obtain 
\begin{align*}
	s_{1}(n,1)&=s_{1}(n,1)\circ 0\\
	          &=1\cdot(D^{n}(\log)\circ 1)=(-1)^{n-1}(n-1)!. 
\end{align*}	
\end{exms}

I use the term \emph{multivariate Stirling polynomial} (MSP) to denote both $\spol{n}{k}$ (MSP of the first kind) and $\bpol{n}{k}$ (MSP of the second kind). \proposition{prop3.7} may be regarded as a good reason for this eponymy (see also my comments in Section~1 concerning notation and terminology).
\begin{rem}\label{rem3.6}
The $\bpol{n}{k}$ are widely known as \emph{partial} exponential Bell polynomials (see e.\,g. \cite{char2002,comt1974}, also the explicit formula \eqref{equation1.4}). Their \emph{complete} form is defined by $B_{n}:=\sum_{k=1}^{n}\bpol{n}{k}$. Applying \proposition{prop3.5}\,(iii) to each term of this sum gives the differential recurrence $B_{n+1}=X_{1}B_{n}+\sum_{j=1}^{n}X_{j+1}\pder{B_{n}}{X_{j}}$, which originally has been studied by Bell \cite{bell1934}. In \cite[p.\,49]{rior1958} the complete Bell polynomials are tabulated up to $n=8$. 
\end{rem}
\raggedbottom
\vspace*{1ex}
\begin{rem}\label{rem3.7}
In \cite[p.\,52 and pp.\,481--483]{knut1997} Knuth analyzes, from a combinatorial point of view, the coefficients of $\bpol{n}{k}$ in connection with \eq{equation3.1} thus establishing a recurrence for $\cbe{n}{k}$. Set $\cbe{n}{k}(\ldots)=0$ for \pt s $\notin\ptsi{n}{k}$. Then $\cbe{1}{1}(1)=1$, and for every $(r_{1},\dots,r_{n-k+2})\in\ptsi{n+1}{k}$: 

\vspace*{-2ex}
\begin{align*}
\cbe{n+1}{k}(r_{1},\dots,r_{n-k+2})&=\cbe{n}{k-1}(r_{1}-1,r_{2},\dots,r_{n-k+2})+\\
	&~~~~\sum_{j=1}^{n-k+1}(r_{j}+1)\cbe{n}{k}(\dots,r_{j}+1,r_{j+1}-1,\dots).
\end{align*}

\vspace*{-1ex}
\sloppy
This could also be obtained more formally by combining \proposition{prop3.5}\,(iii) with \corollary{cor3.6}\,(ii). 

\fussy
I give here without proof also a recurrence for $\cst{n}{k}$ (to be obtained by using \proposition{prop3.5}\,(ii) and \corollary{cor3.6}\,(i)). If we agree in an analogous way to let $\cst{n}{k}$ vanish for \pt s $\notin\ptsi{2n-1-k}{n-1}$, then we have $\cst{1}{1}(1)=1$, and for every $(r_{1},\dots,r_{n-k+2})\in\ptsi{2n+1-k}{n}$:
\begin{alignat*}{2}
\cst{n+1}{k}(&r_{1},\dots,r_{n-k+2})~=~\cst{n}{k-1}(r_{1}-1,r_{2},\dots,r_{n-k+2})\\
                            &+\,(r_{1}-2n+1)\cst{n}{k}(r_{1},r_{2}-1,\dots,r_{n-k+1})\\
														&+\,\sum_{j=2}^{n-k+1}(r_{j}+1)\cst{n}{k}(r_{1}-1,\dots,r_{j}+1,r_{j+1}-1,\dots,r_{n-k+1}).
\end{alignat*}
\end{rem}
%
% End of text
%%%%%%%%%%%%%%%%%%%%%%%%%%%%
%% Section 4: A brief summary on Bell polynomials
%% Last updated: 2021-01-26
%
\section{A brief summary on Bell polynomials}
\flushbottom
\noindent
Replacing the coefficient $a'_{n,k}\circ\varphi$ in \eq{equation3.1} (cf. the proof of \proposition{prop3.3}) by $\dbell{\varphi}{n}{k}$, we obtain for $n\geq0$
\begin{equation}\label{equation4.1}
	D^n(f\circ\varphi) = \sum_{k=0}^{n}\dbell{\varphi}{n}{k}\cdot(D^{k}(f)\circ\varphi).
\end{equation}
This is, in function-algebraic notation, the well-known \emph{Fa\`{a} di Bruno formula} \eqref{equation1.3} for the higher derivatives of a composite function. Though it has been known for a long time, it may, from a systematic point of view, appear appropriate to briefly examine here some of the related classical results on $\bpol{n}{k}$ within our framework.
\label{McKiernan's chain rule}

%\noindent
Let $F\in\funs[X]$. We denote by $\left[F(\varphi)\mid\varphi=0\,\right]$ the result of substituting 0 for $\varphi$ in the monomial products $\varphi^j$ of $F(\varphi)$ with $j\geq1$. Example: Let $F=\varphi+X^2+3(1-X)D(\varphi)$; then $\left[F(\varphi)\mid\varphi=0\,\right]=3D(\varphi)$. 

\begin{prop}\label{prop4.1}
For $1\leq k\leq n$ we have
\[
\dbell{\varphi}{n}{k} = \frac{1}{k!}\left[D^{n}(\varphi^{k})\mid\varphi=0\,\right].
\]
\end{prop}
\begin{proof}
By \eq{equation4.1} we obtain
{\allowdisplaybreaks
\begin{alignat}{2}
	D^{n}(\varphi^{k}) &= D^{n}(\id^{k}\circ\varphi) = \sum_{j=0}^{n}\dbell{\varphi}{n}{j}\cdot(D^{j}(\id^{k})\circ\varphi)\notag \\ 
	  								 &=  \sum_{j=1}^{k}\dbell{\varphi}{n}{j}\cdot\powfall{k}{j}\,\varphi^{k-j} \notag \\ 
										 &= k!\cdot\dbell{\varphi}{n}{k} + \varphi\sum_{j=1}^{k-1}\dbell{\varphi}{n}{j}\cdot\powfall{k}{j}\,\varphi^{k-1-j}.\label{equation4.2}
\end{alignat}
}% end of \allowdisplaybreaks
Taking $\varphi$ to $0$ gives the asserted.
\end{proof}
\begin{prop}\label{prop4.2}
For $1\leq k\leq n$ we have
\[
	\dbell{\varphi}{n}{k}=\frac{1}{k!}\sum_{j=1}^{k}(-1)^{k-j}\binom{k}{j}\varphi^{k-j}D^{n}(\varphi^{j}).
\]
\end{prop}
\begin{proof}
We rewrite \eq{equation4.2} in the form
\begin{align*}
	D^{n}(\varphi^{k})&=\sum_{j=1}^{k} \binom{k}{j}j!\dbell{\varphi}{n}{j}\,\varphi^{k-j}.
\intertext{Applying binomial inversion (cf. \cite[p.\,96-97]{aign1979}) to the equivalent equation}
	\varphi^{-k}D^{n}(\varphi^{k})&=\sum_{j=1}^{k} \binom{k}{j}\varphi^{-j}j!\dbell{\varphi}{n}{j}
\intertext{yields}
		\varphi^{-k}k!\dbell{\varphi}{n}{k}&=\sum_{j=1}^{k}(-1)^{k-j}\binom{k}{j}\varphi^{-j}D^{n}(\varphi^{j}).
\end{align*}
Division by $\varphi^{-k}k!$ then gives the asserted.
\end{proof}

\noindent
Now recall the definition of the subset function $\bru$ in \eq{bruno}.

\begin{prop}\label{prop4.3}
We have $\cbe{n}{k}(r_{1},\ldots,r_{n-k+1})=\bru(r_{1},\dots,r_{n-k+1})$ for all $(r_{1},\dots,r_{n-k+1})\in\ptsi{n}{k}$, that is,
\[
\bpol{n}{k}=\mspace{-6mu}\sum_{\ptsi{n}{k}}\mspace{-5mu}\frac{n!}{r_{1}!\cdots r_{n-k+1}!\cdot1!^{\,r_{1}}\cdots (n-k+1)!^{\,r_{n-k+1}}}X_1^{r_1}\cdots X_{n-k+1}^{r_{n-k+1}}.
\]
\end{prop}
\begin{proof}
It follows from the Leibniz rule \eqref{equation2.4}
\begin{equation*}
D^n(\varphi^{k})=\sum_{\substack{j_{1}+\dots+j_{k}=n \\ 
			                           j_{1},\dots,j_{k}\,\geq\,0}}\,
	\frac{n!}{j_{1}!\cdots j_{k}!}\;D^{j_{1}}(\varphi)\cdots D^{j_{k}}(\varphi).
\end{equation*}
From this we get by \proposition{prop4.1}
\begin{equation}\label{equation4.3}
\bpol{n}{k}=\frac{1}{k!}\sum_{\substack{j_{1}+\cdots+j_{k}=n \\ 
			                           j_{1},\,\ldots,j_{k}\,\geq\,1}}\,
	\frac{n!}{j_{1}!\cdots j_{k}!}\;X_{j_{1}}\cdots X_{j_{k}}.
\end{equation}
Denote by $r_{m}$ the number of $j$'s equal to $m\in\{1,\dots,n-k+1\}$. Then, each sequence $(j_1,\ldots,j_k)$ in (\ref{equation4.3}) is obtained from its corresponding linearly ordered $k$-tuple \mbox{$i_{1}\leq\dots\leq i_{k}$} by $\frac{k!}{r_1!\cdots r_{n-k+1}!}$ permutations. Hence (\ref{equation4.3}) becomes
\begin{equation*}
\bpol{n}{k}=\sum_{\substack{i_{1}+\cdots+i_{k}=n \\ 
			                           1\leq i_{1}\leq\ldots\leq i_{k}}}\,
	\frac{1}{r_1!\cdots r_{n-k+1}!}\cdot\frac{n!}{i_{1}!\cdots i_{k}!}\;X_{i_{1}}\cdots X_{i_{k}},
\end{equation*}
where $i_{1}!i_{2}!\cdots i_{k}!=1!^{r_{1}}\cdot 2!^{r_{2}}\cdots (n+k+1)!^{r_{n-k+1}}$. This yields the asserted equation, and $\cbe{n}{k}$ (the coefficient function of $\bpol{n}{k}$ according to \corollary{cor3.6}) is shown to agree with $\bru$ on $\ptsi{n}{k}$.
\end{proof}
\label{Table: Partial Bell polynomials} 
\label{A consequence from homogeneity}

\begin{rem}\label{rem4.1}
Some historical comments related to Fa\`{a} di Bruno's formula are given in \cite{john2002}. One example that `deserves to be better known' (Johnson), is a formula stated by G. Scott (1861) (cf. \cite{scot1861} and \cite[p.\,220]{john2002}). \proposition{prop4.1} reformulates it in function-algebraic terms. According to \cite{todo1981}, the expression for $\dbell{\varphi}{n}{k}$ given in \proposition{prop4.2} is due to J.~Ber\-trand \cite[p.\,140]{bert1864}). Instead of `a not so easy induction' (Todorov), its verification needs merely applying binomial inversion to Scott's formula. Finally, taking $\varphi=\exp$ makes Bertrand's formula into a well-known explicit expression for the Stirling numbers of the second kind (cf. \cite[p.\,97]{aign1979}): $s_{2}(n,k)=\frac{1}{k!}\sum_{j=1}^{k}(-1)^{k-j}\binom{k}{j}j^{n}.$
\end{rem}

\vspace*{1ex}
Three corollaries will be useful in later sections.\\[6pt]
\indent Because of the relatively simple structure of the $\bpol{n}{k}$ the partial derivatives in \proposition{prop3.5}\,(iii) can be given a closed non-differential form.

\begin{cor}\label{cor4.4}
\[
	\pder{\bpol{n}{k}}{X_{j}}=\binom{n}{j}\bpol{n-j}{k-1}\qquad(1\leq j\leq n-k+1).
\]
\end{cor}
\begin{proof}
The assertion follows by applying $\pdop{X_{j}}$ to the explicit expression of $\bpol{n}{k}$ in \proposition{prop4.3}. Observe that $\pdop{X_{j}}$ takes each $(n,k)$-\pt\ into a $(n-j,k-1)$-\pt. The details are left to the reader. 
\end{proof}

\begin{cor}\label{cor4.5}
\[
	\bpol{n}{k} = \sum_{r=0}^{k}\binom{n}{r}X_{1}^{r}\bpol{n-r}{k-r}(0,X_{2},\dots,X_{n-k+1}).
\]
\end{cor}
\begin{proof}
Immediate from \proposition{prop4.3}. See also Comtet \cite[p.\,136]{comt1974}.
\end{proof}

\begin{ntn}
(i) The right-hand side of the equation in \corollary{cor4.5} gives rise to the notation $\widetilde{B}_{n,k}:=\bpol{n}{k}(0,X_{2},\dots,X_{n-k+1})$. We call $\widetilde{B}_{n,k}$ \textit{associated Bell polynomial} (or, \textit{associated} MSP of the second kind). The coefficients of $\comt{n}{k}$ count only partitions with no singleton blocks. Note that $\assn{n}{k}:=\comt{n}{k}(1,\dots,1)$ are the well-known associated Stirling numbers of the second kind  \cite{howa1980,rior1958}.
\label{Table: Associated Bell polynomials}

%\sloppy
(ii) We call \textit{(unsigned) Lah polynomial} the expression
\[
L^{+}_{n,k}:=P_{\ord}=\sum_{\ptsi{n}{k}}\ord(r_1,\ldots,r_{n-k+1})\,X_{1}^{r_1}\ldots X_{n-k+1}^{r_{n-k+1}},
\]
where $\ord$ is the order function in \eq{lah}. We have 
\[
	L^{+}_{n,k}(1,\ldots,1)=\frac{n!}{k!}\binom{n-1}{k-1}=:l^{+}(n,k). 
\]
Let $l(n,k)$ denote the signed Lah numbers $(-1)^{n}l^{+}(n,k)$. Then 
\[
	L^{+}_{n,k}((-1)^{1},(-1)^{2},\ldots,(-1)^{n-k+1})=l(n,k), 
\]
which follows from the observation that $r_{1}+r_{3}+r_{5}+\ldots\equiv n\,(\mspace{-12mu}\mod 2)$ holds whenever $(r_{1},\ldots,r_{n-k+1})\in\ptsi{n}{k}$.
\end{ntn}
\label{Unsigned Lah numbers}
\label{Table: Unsigned Lah polynomials}

\fussy
\begin{cor}\label{cor4.6}
$L^{+}_{n,k}=\bpol{n}{k}(1!X_{1},2!X_{2},\ldots,(n-k+1)!X_{n-k+1})$. 
\end{cor}

\begin{proof}
Immediate from \proposition{prop4.3}. See also Comtet \cite[p.\,134]{comt1974}.
\end{proof}
%
% End of text
%%%%%%%%%%%%%%%%%%%%%%%%%%%%
%% Section 5: Inversion formulas and recurrences
%% Last updated: 2021-01-26
%%
\section{Inversion formulas and recurrences}
\noindent
We now establish some statements concerning inversion as well as recurrence relations of $\sfun{n}{k}$ and $\bpol{n}{k}$. The first one is the polynomial analogue of the well-known inversion law satisfied by the Stirling numbers of the first and second kind (see \eq{equation1.9}).

\begin{thm}[Inversion Law]\label{thm5.1}
For all $n\geq k\geq 1$
\[
\quad\sum\limits_{j=k}^{n} \sfun{n}{j}\bpol{j}{k}=\kron{n}{k} \qquad\text{and}\qquad\sum\limits_{j=k}^{n} \bpol{n}{j}\sfun{j}{k}=\kron{n}{k}.
\]
\end{thm}

\sloppy
\begin{rem}\label{rem5.1}
Defining lower triangular matrices $\mata_{n}:=(\sfun{i}{j})_{1\leq i,j\leq n}$ and $\matb_{n}:=(\bpol{i}{j})_{1\leq i,j\leq n}$ we can rewrite the statements of the theorem as matrix inversion formulas, for instance, the first one: $\mata_{n}\matb_{n}=\mati_{n}$ (identity matrix) for every $n\geq 1$ (which may also be equivalently expressed by means of differential terms: $\mata_{n}^{\varphi}\,\matb_{n}^{\varphi}=\mati_{n}$). In the special case $X_{1}=X_{2}=\ldots=1$, where the Stirling numbers of the first and second kind are the entries of $\mata_{n}$ and $\matb_{n}$, respectively, both matrices can be considered as transformation matrices connecting the linearly independent polynomial sequences $(x^{1},\ldots,x^{n})$ and $(\powfall{x}{1},\ldots,\powfall{x}{n})$ (cf. \cite{aign1979}). Unfortunately, there is much to suggest that this method does not work in our general case.  The following proof therefore makes no use of it.
\end{rem}

\fussy
\begin{proof}
We prove the first equation of \theorem{thm5.1}. Suppose $1\leq k\leq n$ and denote by $d_{n,k}$ the sum $\sum_{j=k}^{n}a_{n,j}b_{j,k}$. We need to show that $d_{n,k}=\kron{n}{k}$. This is clear for $n=1$. Now we proceed by induction on $n$ using the differential recurrences in the Propositions \ref{prop3.1} and \ref{prop3.3}. First, observe that applying $D$ to both sides of the induction hypothesis yields $D(d_{n,k})=D(\kron{n}{k})=0$, whence
\begin{equation}\label{equation5.1}
	\sum_{j=k}^{n} D(a_{n,j})b_{j,k} = -\sum_{j=k}^{n} a_{n,j}D(b_{j,k}).
\end{equation}
We then have
{\allowdisplaybreaks
\begin{align*}
	d_{n+1,k} &= a_{n+1,n+1}\,b_{n+1,k} + \sum_{j=k}^{n} a_{n+1,j}\,b_{j,k}\\
	          &= \frac{a_{n,n}}{D(\varphi)}\,b_{n+1,k} +
						     \frac{1}{D(\varphi)}\sum_{j=k}^{n}(a_{n,j-1}\,b_{j,k}
								 +D(a_{n,j})b_{j,k})             													 \\
					  &= a_{n,n}\left(b_{n,k-1}+\frac{D(b_{n,k})}{D(\varphi)}\right) +
						     \frac{1}{D(\varphi)}\sum_{j=k}^{n}(a_{n,j-1}D(\varphi)b_{j-1,k-1}\,+
																																					 \\	
						&\quad\qquad\qquad\qquad + a_{n,j-1}D(b_{j-1,k})) + \frac{1}{D(\varphi)}\sum_{j=k}^{n}D(a_{n,j})b_{j,k}.
\end{align*}
}% end of allowdisplaybreaks

\noindent
Replacing the last sum by the right-hand side of \eq{equation5.1}, we obtain after a short computation:
\begin{alignat*}{2}
	d_{n+1,k} &= a_{n,n}\,b_{n,k-1}+ a_{n,n}\frac{D(b_{n,k})}{D(\varphi)}+\sum_{j=k}^{n} a_{n,j-1}\,b_{j-1,k-1}- a_{n,n}\frac{D(b_{n,k})}{D(\varphi)}\\
						&= \sum_{j=k-1}^{n} a_{n,j}\,b_{j,k-1}= \kron{n}{(k-1)} = \kron{(n+1)}{k}.\hspace*{9.6em}\qedhere
\end{alignat*}
\end{proof}

We conclude from \theorem{thm5.1} a statement that generalizes `Stirling inversion' for sequences of real numbers; cf. \cite[Corollary 3.38\,(ii)]{aign1979}.

\sloppy
\begin{cor}[Inversion of sequences]\label{cor5.2}
Let $\mathcal{E}$ be an arbitrary overring of $\ints[X_{1},\ldots,X_{n-k+1}]$, $P_{0},P_{1},P_{2},\ldots$ and $Q_{0},Q_{1},Q_{2},\ldots$ any sequences in $\mathcal{E}$. Then the following conditions are equivalent:

\vspace*{-2ex}
\fussy
\begin{alignat*}{3}
	\text{\textup{(i)}}\qquad  P_{n} &= \sum_{k=0}^{n}\bpol{n}{k}Q_{k}&~~~\text{for all~}n\geq0,\\
	\text{\textup{(ii)}}\qquad Q_{n} &= \sum_{k=0}^{n}\sfun{n}{k}P_{k}&~~~\text{for all~}n\geq0.
\end{alignat*}
\end{cor}

\fussy
\begin{exms}\label{ex5.2}
(i) \theorem{thm5.1} implies the above mentioned special inversion law for the Stirling numbers $s_{1}(n,k)=\sfun{n}{k}(1,\ldots,1)$ and $s_{2}(n,k)=\bpol{n}{k}(1,\ldots,1)$ of the first and second kind (see \proposition{prop3.7}): 
\[
	\sum_{j=k}^{n}s_{1}(n,j)s_{2}(j,k)=\kron{n}{k}\qquad(1\leq k\leq n).
\]
\noindent
(ii) The signed Lah numbers are known to be self-inverse:
\[
	\sum_{j=k}^{n}l(n,j)l(j,k)=\kron{n}{k}\qquad(1\leq k\leq n).
\]
In order to prove this, set $\varphi(x)=-\frac{x}{1+x}$. Since $\varphi$ is involutory ($\inv{\varphi}=\varphi$) and $\varphi(0)=0$, \proposition{prop3.4}\,(ii) yields 
\[
	\dsfun{\varphi}{n}{k}(0)=(\dbell{\inv{\varphi}}{n}{k}\circ\varphi)(0)=\dbell{\varphi}{n}{k}\circ\varphi(0)=\dbell{\varphi}{n}{k}(0). 
\]
It is easily seen that $D^{j}(\varphi)(0)=(-1)^{j}j!$ for all $j\geq 1$. Thus by \corollary{cor4.6} we have $\dbell{\varphi}{n}{k}(0)=\bpol{n}{k}(-1!,2!,-3!,4!,\ldots)=l(n,k)$. Applying \theorem{thm5.1} then gives the desired result.
\end{exms}
\label{Signed Lah numbers}

\begin{thm}\label{thm5.3}
For $1\leq k\leq n$
\begin{alignat*}{2}
	\text{\textup{(i)}} \qquad \sfun{n}{k}&= \bpol{n}{k}(\sfun{1}{1},\dots,\sfun{n-k+1}{1}),\\
  \text{\textup{(ii)}}\qquad \bpol{n}{k}&= \sfun{n}{k}(\sfun{1}{1},\dots,\sfun{n-k+1}{1}).
\end{alignat*}	
\end{thm}

\begin{proof}
(i): By \proposition{prop3.4}\,(ii) $\dsfun{\varphi}{n}{k}=\dbell{\inv{\varphi}}{n}{k}\circ\varphi$. \corollary{cor3.6} yields
\[
	\dbell{\inv{\varphi}}{n}{k}\circ\varphi = \sum_{\ptsi{n}{k}}\mspace{8mu}\cbe{n}{k}(r_1,\ldots,r_{n-k+1})\cdot\prod_{j=1}^{n-k+1} (D^{j}(\inv{\varphi})\circ\varphi)^{r_{j}}.
\]
Setting $k=1$ we get $D^{j}(\inv{\varphi})\circ\varphi=\dbell{\inv{\varphi}}{j}{1}\circ\varphi=\dsfun{\varphi}{j}{1}$ for every $j\geq 1$, hence
\begin{align*}
	\dsfun{\varphi}{n}{k} &= \sum_{\ptsi{n}{k}}\mspace{8mu}\cbe{n}{k}(r_1,\ldots,r_{n-k+1})(\dsfun{\varphi}{1}{1})^{r_{1}}\cdots(\dsfun{\varphi}{n-k+1}{1})^{r_{n-k+1}}\\
	                      &=\bpol{n}{k}(\dsfun{\varphi}{1}{1},\dots,\dsfun{\varphi}{n-k+1}{1})\\
												&=\bpol{n}{k}(A_{1,1},\ldots,A_{n-k+1,1})^{\varphi}.
\end{align*}
(ii): Similarly by \proposition{prop3.4}\,(i).
\end{proof}

The equations of \theorem{thm5.3} can be rewritten as statements about the polynomials $S_{n,k}$.

\begin{cor}\label{cor5.4}
\begin{alignat*}{3}
\text{\textup{(i)}} \qquad\spol{n}{k}&=X_{1}^{k-1}\cdot\bpol{n}{k}(\spol{1}{1},\dots,\spol{n-k+1}{1}),                                                           \\
\text{\textup{(ii)}}\qquad\bpol{n}{k}&=X_{1}^{2k-n}\cdot\spol{n}{k}(\spol{1}{1},\dots,\spol{n-k+1}{1}).
\end{alignat*}	
\end{cor}

\begin{proof}
We show only (ii). It follows from \theorem{thm5.3}\,(ii)
{\allowdisplaybreaks
\begin{alignat*}{3}
	\bpol{n}{k} 
		&= \frac{1}{(\sfun{1}{1})^{2n-1}}\,\spol{n}{k}\left(\frac{\spol{1}{1}}{X_{1}^{1}},\ldots,\frac{\spol{n-k+1}{1}}{X_{1}^{2(n-k+1)-1}}\right) = X_{1}^{2n-1}\,\times\\
		&\mspace{-8mu}\sum_{\ptsi{2n-1-k}{n-1}}\mspace{-5mu}\left(\cst{n}{k}(r_1,r_{2},\ldots)\cdot\mspace{-8mu}\prod_{j=1}^{n-k+1}(\spol{j}{1})^{r_{j}}\cdot X_{1}^{-\sum_{j=1}^{n-k+1}(2j-1)r_{j}}\right)\\
		&= X_{1}^{2n-1}\cdot X_{1}^{-(3n-2k-1)}\cdot\,\spol{n}{k}(\spol{1}{1},\ldots,\spol{n-k+1}{1}).
\end{alignat*}
In the last two lines \corollary{cor3.6}\,(ii) has been used.
}% end of allowdisplaybreaks
\end{proof}

\begin{rem}\label{rem5.3}
Though equations (i) from both \theorem{thm5.3} and \corollary{cor5.4} look like recurrences, their practical (computational) value is rather poor, insofar as they work recursively only for $2\leq k\leq n$. For example, one can actually get $\spol{5}{3}$ by evaluating $X_{1}^{2}\cdot\bpol{5}{3}(\spol{1}{1},\spol{2}{1},\spol{3}{1})$ to $45X_{1}^{2}X_{2}^{2}-10X_{1}^{3}X_{3}$. It should be noted, however, that the very first members of each generation, $\spol{n}{1}$, are the most complicated, and in this case, of all things, (i) yields the empty statement $\spol{n}{1}=\bpol{n}{1}(\spol{1}{1},\dots,\spol{n}{1})$, where $\bpol{n}{1}=X_{n}$.
\end{rem}

\begin{rem}\label{rem5.4}
Through \corollary{cor5.4}, the particular role of $X_{1}$ becomes evident. The exponents appearing here have a combinatorial meaning. Given any $(r_{1},\dots,r_{n-k+1})\in\ptsi{2n-1-k}{n-1}$, one has $r_{1}\geq k-1\geq 0$. For \pt s from $\ptsi{n}{k}$, a corresponding, however possibly negative lower bound holds: $r_{1}\geq 2k-n$. 
\end{rem}

\sloppy
\begin{exm}\label{ex5.5}
It has already been illustrated that putting $X_{j}=1$ $(j=1,2,\dots)$ converts MSP relations into statements about Stirling numbers. So we may ask what in this regard \theorem{thm5.3}\,(i) is about. By \proposition{prop3.7} we obtain a neat identity for the signed Stirling numbers of the first kind:

\fussy
%\vspace*{-4ex}
\begin{equation}\label{equation5.2}
	s_{1}(n,k)=\bpol{n}{k}(s_{1}(1,1),\dots,s_{1}(n-k+1,1)).
\end{equation}

\sloppy
\noindent
We know that $\bpol{n}{k}=P_{\bru}$ (\proposition{prop4.3}) and $s_{1}(j,1)=$ \mbox{$(-1)^{j-1}(j-1)!$} (\examples{ex3.5}\,(iii)). Hence a straightforward evaluation of the right-hand side of \eqref{equation5.2} eventually yields $s_{1}(n,k)=(-1)^{n-k}c(n,k)$ with $c(n,k)=\sum_{\ptsi{n}{k}}\cau(r_{1},\dots,r_{n-k+1})$ (Cauchy's famous enumeration of $n$-permutations with exactly $k$ cycles by means of the cycle function $\cau$, (\ref{cauchy})). Compare equation [3i] in \cite[p.\,135]{comt1974} for a signless version of \eq{equation5.2}.
\end{exm}

\fussy
We will now establish \emph{one} recurrence relation that is satisfied by both $\sfun{n}{k}$ and $\bpol{n}{k}$. However, as with \theorem{thm5.3}, the recurrence does not work for $k=1$, since $\sfun{j}{1}$ and $\bpol{j}{1}$ ($1\leq j\leq n-k+1$) are needed as initial values.
%\smallskip

\begin{prop}\label{prop5.5}
Let $n,k$ be integers with $1\leq k\leq n$. Then we have
\begin{alignat*}{2}
	\text{\textup{(i)}} \qquad \sfun{n}{k} = \sum_{j=1}^{n-k+1}\binom{n-1}{j-1}\sfun{j}{1}\sfun{n-j}{k-1},\\
  \text{\textup{(ii)}}\qquad \bpol{n}{k} = \sum_{j=1}^{n-k+1}\binom{n-1}{j-1}\bpol{j}{1}\bpol{n-j}{k-1}.
\end{alignat*}
\end{prop}

\sloppy
\begin{proof}
(i) can be easily inferred from (ii): We transform $B$ into $B^{\inv{\varphi}}$ and apply $\circ\,\varphi$ (from the right) on both sides of the equation. Then, \proposition{prop3.4}\,(ii) yields the desired statement.\,---\,(ii): Eliminating the partial derivative in \proposition{prop3.5}\,(iii) by \corollary{cor4.4} leads to
\begin{align*}
	\bpol{n}{k} &= X_{1}\bpol{n-1}{k-1} + \sum_{j=1}^{n-1-k+1}X_{j+1}\binom{n-1}{j}\bpol{n-1-j}{k-1}\\
	            &= \sum_{j=1}^{n-k+1}X_{j}\binom{n-1}{j-1}\bpol{n-j}{k-1}. 
\end{align*}
The observation $X_{j}=\bpol{j}{1}$ completes the proof.
\end{proof}

\fussy
\begin{cor}\label{cor5.6}
\[
  \comt{n}{k} = \sum_{j=2}^{n-k+1}\binom{n-1}{j-1}X_{j}\comt{n-j}{k-1} \qquad(1\leq k\leq n).
\]
\end{cor}

\begin{rem}\label{rem5.6}
\proposition{prop5.5}\,(ii) is stated in Charalambides \cite{char2002} (together with a generating function proof; see ibid., p.\,415). In the special case $k=2$ we have $\bpol{j}{1}\bpol{n-j}{k-1}=X_{j}X_{n-j}$. Thus (ii) becomes a simple explicit formula for $\bpol{n}{2}$; see also \cite[p.\,48]{rior1958}.
\end{rem}

\sloppy
\begin{rem}\label{rem5.7}
From \corollary{cor5.6} it follows by an inductive argument that \mbox{$\comt{2n-l}{n}=0$}, if $l\geq1$, and $\comt{2n}{n}=(1\cdot 3\cdot 5\cdots(2n-1))X_{2}^{n}$.
\end{rem}

\fussy
\begin{exms}[\proposition{prop5.5}]\label{ex5.8}
The substitution $X_{j}=1$ makes (i) and (ii) into statements about Stirling numbers. Observe $s_{2}(j,1)=1$ and $s_{1}(j,1)=(-1)^{j-1}(j-1)!$ (\examples{ex3.5}\,(iii)). Then, a little calculation yields the following well-known identities (cf. \cite[p.\,42-43]{rior1958} and \cite[p.\,68]{knut1997}):
\begin{equation*}\hfill
\cycs{n+1}{k+1}=\sum\limits_{j=k}^{n}\frac{n!}{j!}\cycs{j}{k}, \qquad\qquad 
	\subs{n+1}{k+1}=\sum\limits_{j=k}^{n}\binom{n}{j}\subs{j}{k}. 
\end{equation*}	
\end{exms}
%%
%% End of text
%%%%%%%%%%%%%%%%%%%%%%%%%%%%
%% Section 6: Explicit formulas for S_{n,k}
% Update: 2021-01-26
%
\section{Explicit formulas for $S_{n,k}$}
\noindent
We now pass to the task of finding a fully explicit expression for $\spol{n}{k}$. At first glance it seems a viable idea to get $\spol{n}{k}$ by elimination from the inversion law (\theorem{thm5.1}), which is a linear system already in triangular form. In fact, one can verify by induction that every MSP of the first kind can be expressed in terms of Bell polynomials. For instance, in the leading case $k=1<n$ we obtain 
\begin{equation}\label{equation6.1}
\begin{aligned}
	\spol{n}{1}=&-X_{1}^{n-2}\bpol{n}{1}+\\
	            &\mspace{-25mu}\sum_{r=1}^{n-2}(-1)^{r+1}\mspace{-25mu}\sum_{1<j_{1}<\dots<j_{r}<n}\mspace{-15mu}X_{1}^{(n-2)-(j_{1}+\dots+j_{r})}\bpol{n}{j_{r}}\bpol{j_{r}}{j_{r-1}}\cdots\bpol{j_{1}}{1}.
\end{aligned}							
\end{equation}
A quite similar version of \eq{equation6.1} (with $X_{1}=1$) together with an evaluation of the case $n=5$ is mentioned by Figueroa and Gracia-Bond\'{i}a in connection with the antipode on a Hopf algebra (cf. \cite[equation~(7.8)]{figr2005}). 

While it seems dubious if and how \eq{equation6.1} could be further simplified so as to become practicable, it clearly underlines at least again the particular role of $X_{1}$ already observed in \remark{rem5.4}. So, it may appear a promising idea trying to expand $\spol{n}{k}$ into a finite series 
\[
C_{n,k,1}X_{1}+C_{n,k,2}X_{1}^{2}+C_{n,k,3}X_{1}^{3}+\cdots, 
\]
whose coefficients $C_{n,k,r}$ do neither contain $X_{1}$ nor products of two or more Bell polynomials. 

The main result we are going to establish in \theorem{thm6.1} may be regarded as a non-trivial counterpart of \corollary{cor4.5}; it indeed expresses \emph{all} $\spol{n}{k}$ (and consequently $A_{n,k}$) in terms of associated Bell polynomials.
\begin{thm}\label{thm6.1}
For $n\geq k\geq 1$
\[
	\spol{n}{k} = \sum_{r=k-1}^{n-1}(-1)^{n-1-r}\binom{2n-2-r}{k-1}X_{1}^{r}\comt{2n-1-k-r}{n-1-r}.
\]
\end{thm}
\begin{proof}
The proof is divided into two parts. First we will show by induction\footnote{The reciprocity law for Stirling polynomials, which is proved in \cite{schr2020}, results in a new proof that is independent of \theorem{thm6.1} and does not use induction.} that there are polynomials $C_{n,k,r}\in\ints[X_{2},\ldots,X_{n-k+1}]$ such that
\begin{equation}\label{equation6.2}
	\spol{n}{k} = \sum_{r=k-1}^{n-1} C_{n,k,r}X_{1}^{r},
\end{equation}
where the $C_{n,k,r}$ are uniquely determined by a certain differential recurrence. Therefore, in the second step, it remains to show that the recurrence is satisfied by the coefficients of $X_{1}^{r}$ in the asserted equation of the theorem.

\sloppy
1. For $k=n$ we have $\spol{n}{n}=X_{1}^{n-1}$ (\remark{rem3.3}), and taking $C_{n,n,n-1}=1$ satisfies \eq{equation6.2}. This includes the case $n=1$. Now let $n\geq 1$ and suppose \eq{equation6.2} holds for all $k\in\left\{1,\ldots,n\right\}$, where $C_{n,k,r}\in\ints[X_{2},\ldots,X_{n-k+1}]$, $k-1\leq r\leq n-1$. \proposition{prop3.5}\,(ii) yields $\spol{n+1}{k}=T_{n,k}^{(a)}-T_{n,k}^{(b)}+T_{n,k}^{(c)}$ with
{\allowdisplaybreaks
\begin{align*}
	T_{n,k}^{(a)}&=X_{1}\spol{n}{k-1},\quad	T_{n,k}^{(b)}=(2n-1)X_{2}\spol{n}{k}, \text{~~and~}\\
	T_{n,k}^{(c)}&=\sum_{j=1}^{n-k+1}X_{1}X_{j+1}\pder{\spol{n}{k}}{X_{j}}.
%\end{align*}
%%
\intertext{By the induction hypothesis}
%%
%\begin{align*}
	T_{n,k}^{(a)} &= X_{1}\sum_{r=k-2}^{n-1}C_{n,k-1,r}X_{1}^{r}=\sum_{r=k-1}^{n}C_{n,k-1,r-1}X_{1}^{r},\\
	T_{n,k}^{(b)} &= \sum_{r=k-1}^{n-1}(2n-1)X_{2}C_{n,k,r}X_{1}^{r},\\
	T_{n,k}^{(c)} &= X_{1}X_{2}\sum_{r=k-1}^{n-1}\frac{\partial}{\partial X_{1}}(C_{n,k,r}X_{1}^{r})\\
	&~~~~\qquad+ \sum_{j=2}^{n-k+1}\sum_{r=k-1}^{n-1}X_{1}X_{j+1}\frac{\partial}{\partial X_{j}}(C_{n,k,r}X_{1}^{r})\\
	              &= \sum_{r=k-1}^{n-1}r X_{2}C_{n,k,r}X_{1}^{r} +
								\sum_{r=k}^{n}\left(\sum_{j=2}^{n-k+1}X_{j+1}\frac{\partial C_{n,k,r-1}}{\partial X_{j}}\right)X_{1}^{r}.
\end{align*}
} % End allowdisplaybreaks

\noindent
For the desired expansion of $\spol{n+1}{k}$, we now have to find polynomials $C_{n+1,k,r}\in\ints[X_{2},\ldots,X_{n-k+2}]$  that satisfy the recurrence
\vspace*{-1.2ex}
\begin{equation}\label{equation6.3}
	C_{n+1,k,r}=C_{n,k-1,r-1}-(2n\mspace{-1mu}-\mspace{-1mu}1\mspace{-1mu}-\mspace{-1mu}r)X_{2}C_{n,k,r}+\mspace{-8mu}\sum_{j=2}^{n-k+1}\mspace{-5mu}X_{j+1}\frac{\partial C_{n,k,r-1}}{\partial X_{j}}
\end{equation}

\vspace*{-0.1ex}
\fussy\noindent
together with $C_{n,k,r}=1$, if $n=k=r+1$, and $C_{n,k,r}=0$, if $k=0$ or $k\leq r=n-1$ or $r<k-1$. The polynomials $C_{n,k,r}$ are uniquely determined.

2. We set
\begin{equation}\label{equation6.4}
	C_{n,k,r} := (-1)^{n-1-r}\binom{2n-2-r}{k-1}\comt{2n-1-k-r}{n-1-r}.
\end{equation}
It suffices to show that $C_{n,k,r}$ meets the above conditions. First we check the initial values. Since $\comt{2n-1-k-r}{n-1-r}=\comt{2(n-1-r)-l}{n-1-r}$ where $l=k-1-r$, we get by \remark{rem5.7}: $C_{n,k,r}=0$ for $l\geq 1$, that is, for $r<k-1$. The remaining cases are fairly clear.

Next we substitute \eqref{equation6.4} into \eq{equation6.3}. We start with the last summand on the right-hand side of \eq{equation6.3} by first evaluating the partial derivatives:
\[
	\frac{\partial C_{n,k,r-1}}{\partial X_{j}}=(-1)^{n-r}\binom{2n-1-r}{k-1}\frac{\partial}{\partial X_{j}}\comt{2n-k-r}{n-r}.
\]
\corollary{cor4.4} now yields for $j\geq 2$
\[
	\frac{\partial}{\partial X_{j}}\comt{2n-k-r}{n-r}=\binom{2n-k-r}{j}\comt{2n-k-r-j}{n-r-1}.
\]
\sloppy
Substituting the remaining $C$-terms into \eq{equation6.3}, we obtain after some straightforward calculation, in particular cancelling $(-1)^{n-r}$,
{\allowdisplaybreaks
\begin{align}\label{equation6.5}
&\left\{\binom{2n-r}{k-1}-\binom{2n-1-r}{k-2}\right\}\comt{2n+1-k-r}{n-r} =\notag \\
&\qquad\quad\quad(2n-1-r)\binom{2n-2-r}{k-1}X_{2}\comt{2n-1-k-r}{n-1-r}\notag\\
&\qquad\quad+\binom{2n-1-r}{k-1}\sum_{j=2}^{n-k+1}\binom{2n-k-r}{j}X_{j+1}\comt{2n-k-r-j}{n-1-r}.
\end{align}
} % END allowdisplaybreaks
\fussy
Using the identity
\[
	(2n-1-r)\binom{2n-2-r}{k-1}=(2n-k-r)\binom{2n-1-r}{k-1},
\]
equation \eqref{equation6.5} is equivalently reduced to
\[
	\comt{2n+1-k-r}{n-r}=\sum_{j=1}^{n-k+1}\binom{2n-k-r}{j}X_{j+1}\comt{2n-k-r-j}{n-1-r},
\]
which is the statement of \corollary{cor5.6}. This completes the proof.
\end{proof}

\begin{rem}\label{rem6.1}
(i) Taking $k=1$ on the right-hand side of the equation in \theorem{thm6.1}, we obtain Comtet's expansion \eqref{equation1.7} (see \cite[p.\,151]{comt1974}). It should, however, be pointed out that \eq{equation1.7} was intended for a solution of the Lagrange inversion problem (see Section~7 below). In this context, the idea of a connection between $\inv{f}_{n}$ in \eq{equation1.7} and the Lie derivatives $(\theta D)^n$ studied by the same author \cite{comt1973} did not come into play. 

(ii) Sylvester's note \textit{On Reciprocants} (1886) deals with the task `to express the successive derivatives of $x$ in regard to $y$ in terms of those of $y$ in regard to $x$'. The solution is given in the form of the polynomial terms $\spol{n}{1}$, $1\leq n\leq 7$, probably one of the earliest times in the mathematical literature (cf.\,\cite[p.\,38]{sylv1886} and \cite{sylv1854}).
\end{rem}

%\noindent
\theorem{thm6.1} enables us to find closed expressions for the coefficients $\cst{n}{k}$ in \corollary{cor3.6}\,(i), thus leading to an explicit formula for $\spol{n}{k}$, which corresponds to that for the Bell polynomials in \proposition{prop4.3}.
\begin{prop}\label{prop6.2}
For every $(r_{1},\ldots,r_{n-k+1})\in\ptsi{2n-1-k}{n-1}$ we have
\begin{align*}
	\cst{n}{k}(r_{1},\ldots,r_{n-k+1}) &= (-1)^{n-1-r_{1}}\binom{2n-2-r_{1}}{k-1}\,\times\\
	&\qquad\qquad\qquad\times\cbe{2n-1-k-r_{1}}{n-1-r_{1}}(0,r_{2},\ldots)\\
	                               &= (-1)^{n-1-r_{1}}\bru(k-1,r_{2},\ldots,r_{n-k+1}).
\end{align*}
\end{prop}
\label{Table: Stirling polynomials}
\begin{proof}\sloppy
We first consider the coefficients of the associated Bell polynomials. $\comt{2n-1-k-r}{n-1-r}$ is a sum over all \pt s from \mbox{$\ptsi{2n-1-k-r}{n-1-r}$}, whose number $r_1$ of one-element blocks is zero. Let us abbreviate to $\ptas$ the set of these `associated' \pt s. Then, according to \proposition{prop4.3} we can write
\begin{align*}
	\comt{2n-1-k-r}{n-1-r}=\sum_{\ptas}\cbe{2n-1-k-r}{n-1-r}(0,r_{2},&\ldots,r_{n-k+1})\times\\
	&\times\,X_{2}^{r_{2}}\ldots X_{n-k+1}^{r_{n-k+1}}.	                      
\end{align*}												
Observe now that any partition type $(r_{1},\ldots,r_{n-k+1})\in\ptsi{2n-1-k}{n-1}$ can be obtained by putting some $r=r_1\in\left\{k-1,\ldots,n-1\right\}$ in the first place of $(0,r_{2},\ldots,r_{n-k+1})\in\ptas$ (cf. \remark{rem5.4}). This means that the double summation $\sum_{r=k-1}^{n-1}\sum_{\ptas}$ in \theorem{thm6.1} can be replaced by $\sum_{\ptsi{2n-1-k}{n-1}}$. Comparing the coefficients of $X_{1}^{r_{1}}\cdot\ldots\cdot X_{n-k+1}^{r_{n-k+1}}$ in the resulting expression for $\spol{n}{k}$ and in that of \corollary{cor3.6}\,(i) then finally yields the desired $\cst{n}{k}$.
\end{proof}

\fussy
%\noindent
We now add to the list (\ref{lah}--\ref{bruno}) in Section~1 a new mapping $\sti: \ptsa\longrightarrow\ints$, that might be called \textit{Stirling function}, defined by
\begin{equation}\label{equation6.6}
	\sti(r_{1},r_2,r_3,\ldots) = (-1)^{n-1-r_{1}}\frac{(2n-2-r_{1})!}{(k-1)!r_{2}!r_{3}!\cdot\ldots\cdot (2!)^{r_{2}}(3!)^{r_{3}}\cdot\ldots}.
\end{equation}

\sloppy
%\noindent
\proposition{prop6.2} shows that $\cst{n}{k}$ and $\sti$ agree on \mbox{$\ptsi{2n-1-k}{n-1}$}.\,---\,Another succint expression of the relationship between $\bru$ and $\sti$ is the following identity, which holds for all $(2n-1-k,n-1)$-\pt s.
\begin{cor}\label{cor6.3}
\[
	\binom{2n-1-k}{r_{1}}\sti(r_{1},r_{2},\ldots)=(-1)^{n-1-r_{1}}\binom{2n-2-r_{1}}{k-1}\bru(r_{1},r_{2},\ldots).
\]
\end{cor}

\fussy
%\noindent
One could equivalently rewrite this equation as well for $(n,k)$-\pt s as follows:
\[
	\binom{n}{r_{1}}\sti(r_{1},r_{2},\ldots)=(-1)^{k-r_{1}}\binom{2k-r_{1}}{2k-n}\bru(r_{1},r_{2},\ldots).
\]
Here, however, the assumption $2k\geq n$ must be satisfied (see \remark{rem5.4}).
\begin{exm}\label{ex6.2}
The explicit form of $\cst{n}{k}$ can be used to express $s_{2}(n,k)$ as a sum over $\ptsi{2n-1-k}{n-1}$. Suppose $\funs$ is an extended function algebra. Then, from (2') in \examples{ex3.5}\,(i) we obtain $s_{2}(n,k)=\dsfun{\log}{n}{k}\circ 1$. \corollary{cor3.6}\,(i) implies
\[
	\dsfun{\log}{n}{k}=D(\log)^{-(2n-1)}\sum\cst{n}{k}(r_1,r_2,\ldots)D^{1}(\log)^{r_1}D^2(\log)^{r_2}\cdot\ldots
\]
\sloppy
On the right-hand side we have $D^j(\log)=(-1)^{j-1}(j-1)!\cdot\idty^{-j}$ for \mbox{$1\leq j\leq n-k+1$}, and a straightforward calculation using \eq{equation6.6} yields

\fussy
\begin{equation}\label{equation6.7}
	s_2(n,k)=\negthickspace\sum_{\ptsi{2n-1-k}{n-1}}\negthickspace(-1)^{r_1-(k-1)}\frac{\binom{2n-2}{k-1}}{\binom{2n-2}{r_1}}\cdot\cau(r_1,\ldots,r_{n-k+1}).
\end{equation}
\end{exm}

\begin{rem}\label{rem6.3}
It looks like \corollary{cor6.3} expresses a kind of reciprocity law for the coefficient functions $\bru$ and $\sti$. However, it remains open whether there is a connection with the `duality law' for the Stirling numbers\footnote{My paper \cite{schr2020} answers this question in the affirmative.}: $s_{2}(n,k)=c(-k,-n)$ (cf. \cite[p.\,25]{gest1978} and \cite[p.\,412]{knut1992b}). In view of the fact that many MSP relations carry over more or less verbatim to Stirling numbers (see, e.\,g., \remark{rem4.1}, \theorem{thm5.1}, \corollary{cor5.2}, \example{ex5.8}), we could, in the reverse direction, try to establish a reciprocity law for the corresponding polynomials: $\bpol{n}{k}=Z_{-k,-n}$. In the case $n\geq k\geq 0$ we define $Z_{n,k}:=P_{\cau}$ (see (\ref{cauchy}), the sum to be taken over $\ptsi{n}{k}$). The meaning of both $\bpol{n}{k}$ and $C_{n,k}$ can now be extended to arbitrary $n\in\ints$ with the aid of 
\[
\bpol{n}{n-k}=\sum_{j=0}^{k}\binom{n}{k+j}X_{1}^{n-k-j}\comt{k+j}{j}, 
\]
which follows from \corollary{cor4.5} and also applies to the $C_{n,n-k}$ (compare the corresponding statements concerning associative Stirling numbers in \cite{howa1980}; see also \cite[ch.\,8, exercises 10, 12]{char2002}). Again we actually have $s_{2}(n,k)=\bpol{n}{k}(1,\ldots,1)=C_{-k,-n}(1,\ldots,1)=c(-k,-n)$. However, $\bpol{n}{k}=Z_{-k,-n}$ in general turns out to be false. 

\sloppy
Counterexample: $\bpol{5}{2}\neq C_{-2,-5}=105X_{1}^{-8}X_{2}^{3}-120X_{1}^{-7}X_{2}X_{3}+30X_{1}^{-6}X_{4}$, whereas $\bpol{5}{2}(1,1,1,1)=Z_{-2,-5}(1,1,1,1)=15$.
\end{rem}

\fussy
Equation (i) of the theorem to follow is the MSP version of a well-known formula (cf. \example{ex6.4} below) expressing the signed Stirling numbers of the first kind in terms of the Stirling numbers of the second kind.
\begin{thm}[Schl\"omilch type formulas]\label{thm6.4}
For $n\geq k\geq 1$
{\allowdisplaybreaks
\begin{alignat*}{3}
\text{\textup{(i)}}\quad~\sfun{n}{k} &=\mspace{-4mu}\sum_{r=k-1}^{n-1}\mspace{-6mu}(-1)^{n-1-r}\binom{2n-2-r}{k-1}\binom{2n-k}{r+1-k}\,\times\\
			&\qquad\qquad\qquad\qquad\qquad\qquad\times X_{1}^{r-2n+1}\bpol{2n-1-k-r}{n-1-r},\\
\text{\textup{(ii)}}\quad~\bpol{n}{k} &=\mspace{-4mu}\sum_{r=k-1}^{n-1}\mspace{-6mu}(-1)^{n-1-r}\binom{2n-2-r}{k-1}\binom{2n-k}{r+1-k}\,\times\\
			&\qquad\qquad\qquad\qquad\qquad\qquad\times X_{1}^{2n-1-r}\sfun{2n-1-k-r}{n-1-r}.
\end{alignat*}
}% end of \allowdisplaybreaks
\end{thm}
\begin{proof}
For the purpose of formal convenience we adhere to the convention $\bpol{i}{j}=0$ (and consequently $\comt{i}{j}=0$), if $j<0$. Then, the equation of \corollary{cor4.5} can be rewritten in the form
\[
	\bpol{n}{k}=\sum_{j=0}^{n}\binom{n}{j}X_{1}^{n-j}\comt{j}{k-(n-j)},
\]
whence by binomial inversion
\begin{equation}\label{equation6.8}
	\comt{n}{k} =\sum_{j=0}^{k}(-1)^{j}\binom{n}{j}X_{1}^{j}\bpol{n-j}{k-j}.
\end{equation}
The upper limit $n$ has been replaced here by $k$, since we have $\bpol{n-j}{k-j}=0$ for $j>k$. Substituting \eq{equation6.8} into the equation of \theorem{thm6.1} gives
{\allowdisplaybreaks
\begin{align*}
\sfun{n}{k} &= X_{1}^{-(2n-1)}\spol{n}{k}&\qquad\\
            &=\sum_{r=k-1}^{n-1}\sum_{j=0}^{n-1-r}(-1)^{n-1-r-j}\binom{2n-2-r}{k-1}\binom{2n-1-k-r}{j}\times&\\
           &\qquad\qquad\qquad\qquad\qquad\times\,X_{1}^{r-2n+1+j}\bpol{2n-1-k-r-j}{n-1-r-j}.
\intertext{Define now a new index $s=r+j$ so that $k-1\leq s\leq n-1$. Interchanging the order of summation then leads to
}
\sfun{n}{k} &= \sum_{s=k-1}^{n-1}(-1)^{n-1-s}X_{1}^{s-2n+1}\bpol{2n-1-k-s}{n-1-s}~\times\cdots\\
					 &\qquad\cdots\times\underbrace{\sum_{r=k-1}^{s}\binom{2n-2-r}{k-1}\binom{2n-1-k-r}{s-r}}_{(\ast)}.
\intertext{The sum $(\ast)$ can be simplified by using elementary properties of the binomial coefficients. First check that its summand is equal to $\binom{2n-2-r}{2n-2-s}\binom{2n-2-s}{k-1}$. Then, after some calculation, $(\ast)$ is reduced to
}
(\ast)\mspace{-2mu}&=\mspace{-2mu}\binom{2n-2-s}{k-1}\mspace{-8mu}\sum_{r=k-1}^{s}\binom{2n-2-r}{2n-2-s}\mspace{-2mu}=\mspace{-2mu}\binom{2n-2-s}{k-1}\binom{2n-k}{s+1-k}. 
\end{align*}
}% end of allowdisplaybreaks
Finally, renaming the index $s$ to $r$ gives the assertion (i).

\sloppy
Now we derive (ii) from (i). Using \proposition{prop3.4} and \mbox{applying} $P\longmapsto P^{\inv{\varphi}}$ to both sides of (i) makes the left-hand side into \mbox{$\dbell{\varphi}{n}{k}\circ\inv{\varphi}$}, while the term $X_{1}^{r-2n+1}\bpol{2n-1-k-r}{n-1-r}$ on the right becomes $D(\inv{\varphi})^{r-2n+1}\dbell{\inv{\varphi}}{2n-1-k-r}{n-1-r}$. Next, we apply $\circ\varphi$ (from the right) to both sides of (i). It follows
\begin{alignat*}{3}
	\dbell{\varphi}{n}{k} &=\mspace{-4mu}\sum_{r=k-1}^{n-1}\mspace{-6mu}(-1)^{n-1-r}\binom{2n-2-r}{k-1}\binom{2n-k}{r+1-k}~\times\\
	&\qquad\qquad\qquad\quad\,\times(D(\inv{\varphi})^{r-2n+1}\circ\varphi)
	\cdot(\dbell{\inv{\varphi}}{2n-1-k-r}{n-1-r}\circ\varphi).
\end{alignat*}
Observe now that $D(\inv{\varphi})^{r-2n+1}\circ\varphi=D(\varphi)^{2n-1-r}=(X_{1}^{2n-1-r})^{\varphi}$ and, again by \proposition{prop3.4}, $\dbell{\inv{\varphi}}{2n-1-k-r}{n-1-r}\circ\,\varphi=\dsfun{\varphi}{2n-1-k-r}{n-1-r}$. This completes the proof. (The argument amounts to directly applying \theorem{thm5.3}\,(ii), that is, replacing each $X_j$ with $A_{j,1}$ on both sides of (i).)
\end{proof}

\begin{exm}\label{ex6.4}
Specializing \theorem{thm6.4}\,(i) by taking all indeterminates to 1, we immediately obtain \textit{Schl{\"o}milch's formula for the Stirling numbers of the first kind}~\,(see e.\,g. \cite{char2002}):
\begin{align}\label{equation6.9}
	s_{1}(n,k) &= \sum_{r=k-1}^{n-1}(-1)^{n-1-r}\binom{2n-2-r}{k-1}\binom{2n-k}{r+1-k}~\times\qquad\qquad\\
						 &\mspace{135mu}\times\, s_{2}(2n-1-k-r,n-1-r).\notag
	\end{align}

\noindent
Likewise we get from \theorem{thm6.1} a slightly shorter formula of a similar type:
\begin{equation}\label{equation6.10}
	s_{1}(n,k) = \sum_{r=k-1}^{n-1}(-1)^{n-1-r}\binom{2n-2-r}{k-1}\assn{2n-1-k-r}{n-1-r}.
\end{equation}

\vspace*{1ex}
\fussy
\noindent
\eq{equation6.10} runs with smaller numbers than does \eqref{equation6.9}. Compare, for instance, the computation of $s_{1}(7,4)=-84\cdot 90+56\cdot 150 - 35\cdot 45=-735$ by Schl\"omilch's formula \eqref{equation6.9} with that of $s_{1}(7,4)=-84\cdot 15 + 56\cdot 10 - 35\cdot 1= -735$ by \eq{equation6.10}. However, associated Stirling numbers of the second kind do not have quite simple explicit representations (cf. \cite{howa1980}, in particular equation (3.11), ibid.).

Finally it should be noted that the statement (ii) of \theorem{thm6.4} enables the Stirling numbers of the second kind to be represented using the Stirling numbers of the first kind. This result is due to Gould \cite{goul1960}.
\end{exm}
%
% End of text
%%%%%%%%%%%%%%%%%%%%%%%%%%%%
%% Section 7: Remarks on Lagrange inversion
%% Last update: 2021-01-26
%%
\section{Remarks on Lagrange inversion}
Let $\varphi$ be an analytic function (bijective in the real or complex domain), which is given in the form of a power series 
\[
	\varphi(x)=\sum_{n\geq 1}f_{n}\frac{x^n}{n!}
\]
with zero constant term and $f_{1}\neq 0$. Then, the compositional inverse $\inv{\varphi}$ is unique, and one may ask for the coefficients $\inv{f}_{n}$ in the series expansion $\inv{\varphi}(x)=\sum_{n\geq 1}\inv{f}_{n}x^n/n!$. One possible answer to this is Lagrange's famous inversion formula (cf. \cite{comt1974} and \cite{stan1999} for proofs and further references):
\begin{equation}\label{equation7.1}
	\inv{f}_{n} = \left(\frac{d}{dx}\right)^{n-1}\mspace{-7mu}\left(\left(\frac{x}{\varphi(x)}\right)^{n}\right){\biggr|_{x=0}}.
\end{equation}
\sloppy

This innocent looking expression, however, provides in most cases an all but simple method of computing. Many attempts have therefore been made to obtain alternative and more efficient expressions (see, e.\,g., \cite{krat1988}, \cite{todo1981}). Trying to express $\inv{f}_{n}$ as a function of $f_{1},f_{2},\ldots,f_{n}$ is a quite natural approach. Morse and Feshbach \cite[p.\,412]{mofe1953}, for example, employed the Residue Theorem to show that $\inv{f}_{n}$ can be represented as a polynomial expression over all partitions of $n-1$. Comtet \cite{comt1974} derived the remarkable result that the right-hand side of \eq{equation7.1} is equal to $\sum_{k=0}^{n-1}(-1)^k f_1^{-n-k}\bpol{n+k-1}{k}(0,f_2,\ldots,f_n)$. By \theorem{thm6.1} the latter can easily be seen to agree with $\sfun{n}{1}(f_{1},\ldots,f_{n})$, that is, we have for all $n\geq 1$:
\begin{equation}\label{equation7.2}
	\inv{f}_{n} = \sfun{n}{1}(f_{1},\ldots,f_{n}).
\end{equation}
\label{Inversion of power series}

\vspace*{-2ex}
In the following we will deal with some few function-algebraic aspects of Lagrange inversion. It turns out that \eq{equation7.2} can be proved without using \eqref{equation7.1}. We will also see that inverting a function $\varphi$ is nothing else than switching from $\hlieder{\id}{n}{\varphi}(0)$ to $\hlieder{\varphi}{n}{\id}(0)$ in the corresponding series expansions. A~few examples should briefly illustrate the computational aspects.\\[-4pt]

%\noindent 
Let $\funs$ be an extended function algebra, $\varphi\in\funs$. Throughout this section we assume $f_{0}:=\varphi\circ 0=0$, and $f_{1}:=D(\varphi)\circ 0$ to be a unit.\footnote{In the sequel we use the customary notation $g(0)$ instead of $g\circ 0$, $g\in\funs$}  Since the notion of convergence has no place in $\funs$, we make use of formal power series.
\begin{dfn}\label{dfn7.1}
We say that $\varphi$ is an \textit{exponential generating function} of the sequence of constants $f_{0},f_{1},f_{2},\ldots\in\cons$ (symbolically written $\varphi(x)=\sum_{n\geq0}f_{n}x^{n}/n!$), if $D^{n}(\varphi)(0)=f_{n}$ for every $n\geq0$.\footnote{In the case $f_{0}=0$ we take the lower summation limit to 1.}
\end{dfn}

The exponential generating function of the constant sequence $1,1,1,\ldots$ 
\[ 
e^x:=\exp(x):=\sum_{n\geq0}\frac{x^{n}}{n!}
\]
obviously has the properties that one would expect from an exponential (as has been indicated in Section~2.2).

Our basic statement on inversion is now the following
\begin{prop}\label{prop7.1}
\begin{align*}
\text{\textup{(i)}}~\quad \varphi(x) &= \sum_{n\geq 1}\hlieder{\textup{\idty}}{n}{\varphi}(0)\frac{x^{n}}{n!}\, = \,\sum_{n\geq 1}\dbell{\varphi}{n}{1}(0)\frac{x^{n}}{n!},\\
\text{\textup{(ii)}}~\quad \inv{\varphi}(x) &= \sum_{n\geq 1}\hlieder{\varphi}{n}{\textup{\idty}}(0)\frac{x^{n}}{n!}\, =\, \sum_{n\geq 1}\dsfun{\varphi}{n}{1}(0)\frac{x^{n}}{n!}.
\end{align*}
\end{prop}
\begin{proof}
(i): Clearly $D_{\idty}=D$. Hence for all $n\geq 1$ 
\[
	\dbell{\varphi}{n}{1}(0)=(X_{n})^{\varphi}(0)=D^{n}(\varphi)(0)=D^{n}_{\idty}(\varphi)(0).
\]
(ii): By \proposition{prop3.4}\,(ii) we have 
\[
	\dsfun{\varphi}{n}{1}(0)=(\dbell{\inv{\varphi}}{n}{1}\circ \varphi)(0)=(D^{n}(\inv{\varphi})\circ\varphi)(0).
\]
Thus Pourchet's formula yields
\[
	\hlieder{\varphi}{n}{\textup{\idty}}(0)=(D^{n}(\inv{\varphi})\circ\varphi)(0)=D^{n}(\inv{\varphi})(\varphi(0))=D^{n}(\inv{\varphi})(0).\qedhere
\]	
\end{proof}

\vspace*{1ex}
\begin{cor}\label{cor7.2}
For every $n\geq 1$ 
\begin{align*}
\inv{f}_{n}&=\sfun{n}{1}(f_{1},\ldots,f_{n})\\
           &=\frac{1}{f_{1}^{2n-1}}\cdot\mspace{-12mu}\sum_{\ptsi{2n-2}{n-1}}\frac{(-1)^{n-1-r_{1}}\cdot(2n-2-r_{1})!}{r_{2}!\cdots r_{n}!\cdot(2!)^{r_2}\cdots (n!)^{r_n}}\,f_{1}^{r_1}f_{2}^{r_2}\cdots f_{n}^{r_n}.
\end{align*}
\end{cor}
\begin{proof}
\proposition{prop7.1}\,(ii) yields $\inv{f}_{n}=\dsfun{\varphi}{n}{1}(0)=\sfun{n}{1}(f_{1},\ldots,f_{n})$, where $f_{j}=D^{j}(\varphi)(0)$ ($j=1,\ldots,n$). Then, applying \corollary{cor3.6}\,(i) and \eq{equation6.6} for $k=1$ gives the asserted.
\end{proof}

Except for some simple cases, the higher Lie derivatives $\hlieder{\varphi}{n}{\textup{\idty}}$ turn out to be considerably less complex than the Lagrangian terms $D^{n-1}((\id/{\varphi})^{n})$. Nonetheless, the most advantage is presumably to be gained from applying \theorem{thm6.1} (with $k=1$) or \corollary{cor7.2} to the coefficients $f_{j}=D^{j}(\varphi)(0)$ ($j=1,2,\ldots,n$).

Let us consider three examples.

\begin{exm}\label{ex7.1}
Let $\varphi(x)=x e^{-x}$. As is well-known (see, e.\,g., \cite[p.\,23]{stan1999}), the inverse $\inv{\varphi}$ is exponential generating function of the sequence $r(n)$ ($n=1,2,3,\ldots$) of numbers of rooted (labeled) trees on $n$ vertices. The function $\varphi$ seems to be tailored to the application of \eq{equation7.1}, for we have $D^{n-1}((\idty/\varphi)^n)(x)=D^{n-1}(\exp\circ(n\cdot \idty))(x)=n^{n-1}e^{n x}$ and thus readily by \eq{equation7.1}: $r(n)=\inv{f}_{n}=n^{n-1}e^{n x}|_{x=0}=n^{n-1}$.

Compared with this neat calculation, applying $\sfun{n}{1}$ to the coefficients $f_{j}=D^{j}(\varphi)(0)=(-1)^{j-1}j$ ($j=1,\ldots,n$) gives less satisfactory results. \corollary{cor7.2} yields
\enlargethispage{1ex}
\vspace*{-1.2ex}
\begin{equation}\label{equation7.3}
	\sum_{\ptsi{2n-2}{n-1}}(-1)^{r_{1}}\frac{(2n-2-r_{1})!}{r_{2}!r_{3}!\cdots r_{n}!}\cdot\frac{1}{1!^{r_2}2!^{r_3}\cdots(n-1)!^{r_n}}.
\end{equation}
Here a combinatorial argument is needed to see that \eq{equation7.3} equals $n^{n-1}$.
\end{exm}

\begin{exm}\label{ex7.2}
Let $\varphi(x)=e^{x}-1$. While formula \eqref{equation7.1} fails to yield simple expressions, the Taylor series expansion of $\inv{\varphi}=\log\circ(1+\idty)$ is immediately obtained by \proposition{prop7.1}\,(ii). We have $\dsfun{\varphi}{n}{1}(0)=\sfun{n}{1}(1,\ldots,1)=s_1(n,1)=(-1)^{n-1}(n-1)!$ (see \examples{ex3.5}\,(iii)), and hence $\log(1+x)=\sum_{n\geq 1}(-1)^{n-1}x^{n}/n$.
\end{exm}
\vspace*{1ex}

\begin{exm}\label{ex7.3}
Let $\varphi(x)=1+2x-e^{x}$. According to Stanley \cite[p.\,13]{stan1999}, $\inv{\varphi}$ is exponential generating function of the sequence $t(n)$ ($n=1,2,3,\ldots$), where $t(n)$ denotes the number of total partitions of $\{1,\ldots,n\}$; cf. the fourth problem posed by E. Schr{\"o}der (1870) (`arbitrary set bracketings') \cite[p.\,178]{stan1999}. When applied to $\varphi$, the Lagrange formula does not `seem to yield anything interesting' (Stanley). Let us therefore have a look at what can be achieved with the help of the MSPs. 

We have $f_{1}=1, f_{j}=-1 ~(j\geq 2)$, and hence by Corollaries \ref{cor6.3} and \ref{cor7.2}
\[
	t(n)=A_{n,1}(1,-1,\ldots,-1)=\sum_{\ptsi{2n-2}{n-1}}\binom{2n-2}{r_1}^{-1}\bru(r_1,\ldots,r_n).
\]
The same can be expressed as well in terms of associated Stirling numbers of the second kind by applying \theorem{thm6.1}:
\[
	t(n)=\sum_{r=0}^{n-1}\assn{2n-2-r}{n-1-r}.
\]
Regarding the latter formula, the reader is referred to Comtet \cite[p.\,224]{comt1974}.

Alternatively, we can use the fact that $\inv{f}_{n}=\hlieder{\varphi}{n}{\textup{\idty}}(0)$ (cf. \proposition{prop7.1}\,(ii)). This leads to a recursive solution. The repeated application of the Lie derivation $D_{\varphi}$ inductively results in a representation of the form
\begin{align*}\label{equation7.4}
\begin{aligned}
	              ~&\hlieder{\varphi}{n}{\textup{\idty}}(x) = (2-e^x)^{-(2n-1)}T_{n}(x),
	\quad\text{where}\\&T_{1}(x)=1,\quad
	T_{n}(x)=\sum_{k=0}^{n-1}b_{n-1,k}e^{k x}\quad\text{for~}n\geq 2,
\end{aligned}	
\end{align*}
with non-negative integers $b_{n-1,k}$. It follows for all $n\geq 1$
\[
	t(n)=T_{n}(0)=b_{n-1,0}+b_{n-1,1}+\cdots+b_{n-1,n-1}.
\]
Equating now the coefficients of $e^{k x}$ in $\hlieder{\varphi}{n+1}{\idty}(x)=D_{\varphi}(\hlieder{\varphi}{n}{\idty})(x)$ gives the recurrence
\begin{equation*}\label{equation7.5}
\begin{aligned}
		b_{n,k} &= (2n-k)b_{n-1,k-1}+2kb_{n-1,k}\qquad &(1\leq k\leq n),\\
		b_{i,0} &=\kron{i}{0},\quad b_{i,j}=0 &(0\leq i<j).
\end{aligned}	
\end{equation*}
Some special values:
\[
	b_{n,1}=2^{n-1},\quad b_{n,2}=2^{n-1}(2^{n}-n-1),\quad b_{n,n}=n!.
\]
We obtain $t(1)=1$, $t(2)=0+1=1$, $t(3)=0+2+2=4$, $t(4)=0+2^2+2^2(2^3-3-1)+3!=26$.
\end{exm}

\medskip
Our last statement concerns exponential generating functions of the form $\exp\circ(t\cdot\varphi)$, $t\in\cons$.
\begin{prop}\label{prop7.3}
\begin{alignat*}{3}
\text{\textup{(i)}}~\quad e^{t \varphi(x)}&= \sum_{n\geq 0}\left(\sum_{k=0}^{n}\dbell{\varphi}{n}{k}(0)\,t^{k}\right)\frac{x^{n}}{n!},\qquad\\
\text{\textup{(ii)}}~\quad e^{t \inv{\varphi}(x)}&= \sum_{n\geq 0}\left(\sum_{k=0}^{n}\dsfun{\varphi}{n}{k}(0)\,t^{k}\right)\frac{x^{n}}{n!}.
\end{alignat*}
\end{prop}
\begin{proof}
(i): By Fa\`{a} di Bruno's formula \eqref{equation4.1} and \remark{rem2.3}
\[
	D^{n}(e^{t\varphi})=\sum_{k=0}^{n}\dbell{\mspace{2mu}t\varphi}{n}{k}\cdot(D^{k}(\exp)\circ(t\varphi))=e^{t\varphi}\sum_{k=0}^{n}\dbell{\varphi}{n}{k}\,t^{k},
\]
hence
\[
	D^{n}(e^{t\varphi})(0)=e^{t\varphi(0)}\sum_{k=0}^{n}\dbell{\varphi}{n}{k}(0)\cdot(t\circ 0)^{k}.
\]
Since $e^{t\varphi(0)}=e^{t\cdot 0}=e^{0}=1$ and $t\circ 0=t$, we are done.

Note $\inv{\varphi}(0)=0$; then (ii) follows from (i) by \proposition{prop3.4}.
\end{proof}

We conclude with a well-known example.

\begin{exm}\label{ex7.4}
As in \example{ex7.2}, take $\varphi=\exp-1$. From (2) in \examples{ex3.5}\,(i) we have $\dbell{\varphi}{n}{k}(0)=s_{2}(n,k)$. Then by \proposition{prop7.3}\,(i)
\vspace*{-1ex}
\begin{equation*}\label{equation7.6}
	\exp(t(e^x-1))=\sum_{n\geq 0}\left(\sum_{k=0}^{n}s_{2}(n,k)\,t^k\right)\frac{x^n}{n!}.
\end{equation*}
Put $t=1$. This shows that $\exp(e^x-1)$ is exponential generating function of the Bell number sequence $b(n):=\sum_{k=0}^{n}s_{2}(n,k)$ (see \cite[p.\,13]{stan1999}).

We also have $\dsfun{\varphi}{n}{k}(0)=s_{1}(n,k)$ (from (1) in \examples{ex3.5}\,(i)). Since $\inv{\varphi}=\log\circ(1+\id)$, \proposition{prop7.3}\,(ii) yields
\begin{equation}\label{equation7.7}
	\exp(t\,\log(1+x))=\sum_{n\geq 0}\left(\sum_{k=0}^{n}s_{1}(n,k)\,t^k\right)\frac{x^n}{n!}.
\end{equation}
Thus $(1+x)^t$ turns out to be the exponential generating function of the sequence $\sum_{k=0}^{n}s_{1}(n,k)\,t^k$ ($n=1,2,3,\ldots$) (see, e.\,g., \cite[p.\,281]{char2002}).
\end{exm}
%
% End of text
%%%%%%%%%%%%%%%%%%%%%%%%%%%%
%% Section 8: Concluding remarks
%% Last updated: 2020-11-13
%
%\vspace*{2ex}
%
\section{Concluding remarks}
\noindent
The preceding work was primarily intended as an attempt to introduce, within a general function-algebraic setting, the notion of MSP of the first and second kind (the latter being identical to that of partial Bell polynomial). The investigation was focussed on establishing the inverse relationship as well as other fundamental properties of the two polynomial families. The resulting picture shows that the MSPs may be understood as a kind of strong generalizations of the corresponding Stirling numbers. 
\label{Stirling polynomials and Stirling numbers}

Supplementary to this the reader is referred to a package \cite{schr2013} for \textit{Mathematica}\raisebox{5pt}{\scriptsize{\textregistered}} providing functions that generate the expressions $\spol{n}{k}$, $\sfun{n}{k}$, and $B_{n,k}$ together with a substitution mechanism. 

No attempt was made here to develop a combinatorial interpretation of the coefficient function $\sti: \ptsa\longrightarrow\ints$ (in $\spol{n}{k}$ and $A_{n,k}$). Based on \cite{schm1986}, Haiman and Schmitt \cite{hasc1989} have offered a satisfactory explanation of Comtet's expansion \eqref{equation1.7} (cf. \remark{rem6.1}) from an incidence algebra point of view (using colored partitions of finite colored sets). So, it might be worth-while examining whether this idea will also apply to the general case formulated in \theorem{thm6.1}. The `ban on one-element blocks' observed in the case $k=1$ \cite[p.\,180]{hasc1989} may be seen as a sign in this direction, since it is as well a striking feature of the whole family $\spol{n}{k}$ (see, e.\,g., \remark{rem5.4}, \proposition{prop6.2}, and \eq{equation6.6}).
\newpage
% End of text
%%%%%%%%%%%%%%%%%%%%%%%%%%%%
%%
%% ----------------------------------------------------
\clearpage
%% Table MSP
%% Last update: 2013-07-08
%%
\begin{table}
	\begin{center}
	  \renewcommand{\arraystretch}{1.3}
		\begin{tabular}[ht]{| l | l |}
		\hline
		\textbf{MSP of the 1st kind} & \textbf{MSP of the 2nd kind} \\ \hline\hline
		%%-------------------------------------------------------
		$\spol{1}{1}=1$     & $\bell{1}{1}=X_1$   \\ \hline
		%%-------------------------------------------------------
		$\spol{2}{1}=-X_2$  & $\bell{2}{1}=X_2$   \\
		$\spol{2}{2}=X_1$   & $\bell{2}{2}=X_1^2$ \\ \hline
		%%-------------------------------------------------------
		$\spol{3}{1}=3X_{2}^{2}-X_{1}X_{3}$ &  $\bell{3}{1}=X_3$ \\
		$\spol{3}{2}=-3X_{1}X_{2}$		    	&	 $\bell{3}{2}=3X_{1}X_{2}$ \\
		$\spol{3}{3}=X_{1}^{2}$						&	 $\bell{3}{3}=X_{1}^{3}$ \\ \hline
		%%-------------------------------------------------------
		$\spol{4}{1}=-15X_{2}^{3}+10X_{1}X_{2}X_{3}-X_{1}^{2}X_{4}$  &
					$\bell{4}{1}=X_{4}$                                    \\
		$\spol{4}{2}=15X_{1}X_{2}^{2}-4X_{1}^{2}X_{3}$               &
					$\bell{4}{2}=4X_{1}X_{3}+3X_{2}^{2}$                   \\
		$\spol{4}{3}=-6X_{1}^{2}X_{2}$                               &
					$\bell{4}{3}=6X_{1}^{2}X_{2}$                          \\
		$\spol{4}{4}=X_{1}^{3}$                                       &
					$\bell{4}{4}=X_{1}^{4}$                                \\ \hline
		%%-------------------------------------------------------
    $\spol{5}{1}=105X_{2}^{4}-105X_{1}X_{2}^{2}X_{3}+10X_{1}^{2}X_{3}^{2}$ &
           $\bell{5}{1}=X_{5}$                										\\						
		\qquad\quad$+ 15X_{1}^{2}X_{2}X_{4}-X_{1}^{3}X_{5}$         &
					  		 																									\\	
    $\spol{5}{2}=-105X_{1}X_{2}^{3}+60X_{1}^{2}X_{2}X_{3}
		             -5X_{1}^{3}X_{4}$								                &
					$\bell{5}{2}=5X_{1}X_{4}+10X_{2}X_{3}$                  \\
    $\spol{5}{3}=45X_{1}^{2}X_{2}^{2}-10X_{1}^{3}X_{3}$						&
					$\bell{5}{3}=15X_{1}X_{2}^{2}+10X_{1}^{3}X_{3}$         \\
		$\spol{5}{4}=-10X_{1}^{3}X_{2}$                               &
					$\bell{5}{4}=10X_{1}^{3}X_{2}$                          \\
		$\spol{5}{5}=X_{1}^{4}$                                       &
					$\bell{5}{5}=X_{1}^{5}$                                 \\ \hline
		%%-------------------------------------------------------
		$\spol{6}{1}=-945X_{2}^{5}+1260X_{1}X_{2}^{3}X_{3}-280X_{1}^{2}X_{2}X_{3}^{2}$ &
		      $\bell{6}{1}=X_{6}$                                     \\
		\qquad$~\;-210X_{1}^{2}X_{2}^{2}X_{4}+35X_{1}^{3}X_{3}X_{4}$                                                                  &
		                                                              \\ 
		\qquad$~\;+21X_{1}^{3}X_{2}X_{5}-X_{1}^{4}X_{6}$                                     &
			                                                            \\
		$\spol{6}{2}=945X_{1}X_{2}^{4}-840X_{1}^{2}X_{2}^{2}X_{3}+70X_{1}^{3}X_{3}^{2}$																															     &
          $\bell{6}{2}=10X_{3}^{2}+15X_{2}X_{4}+6X_{1}X_{5}$      \\    
		\qquad$~\;+105X_{1}^{3}X_{2}X_{4}-6X_{1}^{4}X_{5}$            & 
		                                                              \\
		$\spol{6}{3}=-420X_{1}^{2}X_{2}^{3}+210X_{1}^{3}X_{2}X_{3}-15X_{1}^{4}X_{4}$ &
		      $\bell{6}{3}=15X_{2}^{3}+60X_{1}X_{2}X_{3}+15X_{1}^{2}X_{4}$ \\
		$\spol{6}{4}=105X_{1}^{3}X_{2}^{2}-20X_{1}^{4}X_{3}$          &
		      $\bell{6}{4}=45X_{1}^{2}X_{2}^{2}+20X_{1}^{3}X_{3}$    \\ 
		$\spol{6}{5}=-15X_{1}^{4}X_{2}$                               &
		      $\bell{6}{5}=15X_{1}^{4}X_{2}$                         \\
		$\spol{6}{6}=X_{1}^{5}$                                       &
					$\bell{6}{6}=X_{1}^{6}$                                \\
		\hline
		\end{tabular}
		\caption{Generations 1\,--\,6 of the multivariate Stirling polynomials of the first kind ($\spol{n}{k}$) and of the second kind ($\bell{n}{k}$: partial exponential Bell polynomials).}\label{table_msp}
	\end{center}
\end{table}	
%%
%% End of table
%%
\appendix
%% Bibliographic references:
\clearpage
%% Bibliographic references
%% msp_schreiber_bib.tex
%% Last update 2013-11-05
%% Revised 2020-12-05
%%-------------------------
%%

%%-----------------------
%%
%%%%%%%%%%%%%%%%%%%%%%%%%%%%%%%%%%%%%%%%%%%%%%%%%%%%%%%%%%%%%%%%%%%%%%%%
%\nocite* %% causes bibtex to take over all bibliographic entries
%%
\end{document}